\documentclass[12pt]{amsart}

\usepackage[toc]{appendix} 
\usepackage{latexsym}
\usepackage{amsfonts}
\usepackage{amssymb}
\usepackage{amsmath}
\usepackage{mathrsfs}
\usepackage{bbm}
\usepackage{dsfont}
\usepackage{hyperref}
\allowdisplaybreaks

\usepackage[toc]{appendix}

\newcommand{\be}{\begin{equation}}
\newcommand{\ee}{\end{equation}}
\newcommand{\bea}{\begin{eqnarray}}
\newcommand{\eea}{\end{eqnarray}}
\newcommand{\bean}{\begin{eqnarray*}}
\newcommand{\eean}{\end{eqnarray*}}
\newcommand{\brray}{\begin{array}}
\newcommand{\erray}{\end{array}}

\newtheorem{dfn}{Definition}[section]
\newtheorem{thm}[dfn]{Theorem}
\newtheorem{lmma}[dfn]{Lemma}
\newtheorem{ppsn}[dfn]{Proposition}
\newtheorem{crlre}[dfn]{Corollary}
\newtheorem{xmpl}[dfn]{Example}
\newtheorem{rmrk}[dfn]{Remark}

\newcommand{\bdfn}{\begin{dfn}\rm}
\newcommand{\bthm}{\begin{thm}}
\newcommand{\blmma}{\begin{lmma}}
\newcommand{\bppsn}{\begin{ppsn}}
\newcommand{\bcrlre}{\begin{crlre}}
\newcommand{\bxmpl}{\begin{xmpl}}
\newcommand{\brmrk}{\begin{rmrk}\rm}

\newcommand{\edfn}{\end{dfn}}
\newcommand{\ethm}{\end{thm}}
\newcommand{\elmma}{\end{lmma}}
\newcommand{\eppsn}{\end{ppsn}}
\newcommand{\ecrlre}{\end{crlre}}
\newcommand{\exmpl}{\end{xmpl}}
\newcommand{\ermrk}{\end{rmrk}}

\newcommand{\bbc}{\mathbb{C}}
\newcommand{\bbz}{\mathbb{Z}}

\newcommand{\bbn}{\mathbb{N}}
\newcommand{\bbr}{\mathbb{R}}

\newcommand{\bbt}{\mathbb{T}}

\newcommand{\clh}{\mathcal{H}}

\newcommand{\clk}{\mathcal{K}}

\newcommand{\clg}{\mathcal{G}}

\makeatletter
\let\@wraptoccontribs\wraptoccontribs
\makeatother

\author{Anbu Arjunan, Sruthymurali and S. Sundar}
\title{KMS STATES ON $C_{c}^{*}(\bbn^2)$}

\date{}
\begin{document}
\maketitle
\begin{abstract}
Let $C_c^{*}(\mathbb{N}^{2})$ be the universal $C^{*}$-algebra generated by a semigroup of isometries  $\{v_{(m,n)}: m,n \in \mathbb{N}\}$ whose range projections commute.  We  analyse the structure of  KMS 
states on $C_{c}^{*}(\bbn^2)$ for the time evolution determined by a homomorphism $c:\mathbb{Z}^{2} \to \mathbb{R}$. In contrast to the reduced version $C_{red}^{*}(\mathbb{N}^{2})$, we show that the set of KMS states on $C_{c}^{*}(\mathbb{N}^{2})$ has a 
rich  structure. In particular, we exhibit uncountably many extremal KMS states of type I, II and III. 
\end{abstract}

 \noindent {\bf AMS Classification No. :} {Primary 46L05; Secondary 46L30, 20M99.}  

\section{Introduction}
In the recent years, there has been a flurry of activity centered around understanding the inner structure (like $K$-theory, KMS states etc...) of $C^{*}$-algebras associated to semigroups. The  revival of the subject of semigroup $C^{*}$-algebras, especially in the last decade,  can be attributed to the 
  works of Cuntz and Li on  $C^{*}$-algebras associated to rings (\cite{Cuntz}, \cite{Cuntz-Li}, \cite{Cuntz-Li-1}).   Ring $C^{*}$-algebras are defined exactly like group $C^{*}$-algebras taking into account both the addition and the multiplication rule of the ring, where the ring involved is usually assumed to be an integral domain. Soon, it was realised that it is only appropriate to view these algebras  as semigroup $C^{*}$-algebras.

Let $P$ be a cancellative semigroup. Denote the left regular representation of $P$ on $\ell^{2}(P)$ by $V:=\{V_{a}\}_{a \in P}$. The reduced $C^{*}$-algebra of the semigroup $P$, denoted $C_{red}^{*}(P)$, is the $C^{*}$-algebra generated by $\{V_{a}:a \in P\}$. Li also defines a universal version, denoted $C^{*}(P)$, generated by isometries $\{v_{a}:a \in P\}$ and projections corresponding to `certain ideals of $P$' satisfying  relations that reflect the relations in the regular representation.  The study of $C_{red}^{*}(P)$  has received much attention in  the recent years. Questions concerning its nuclearity,  $K$-theory and  the structure of KMS states on $C_{red}^{*}(P)$ were investigated intensively and satisfactory answers were obtained for a large class of semigroups. Some of the notable papers that explore the above mentioned issues are \cite{Li-Cuntz-Echterhoff}, \cite{Li13}, \cite{Li-semigroup}, \cite{Cuntz_Laca_KMS}, \cite{Laca_Raeburn_KMS}, \cite{Laca_Neshveyev_KMS}  and \cite{Larsen_Stam}.

However, isometric representations of semigroups other than the regular representation were also considered in the literature and the associated $C^{*}$-algebras were analysed.
Let $G$ be a discrete countable abelian group and let $P \subset G$ be a subsemigroup containing the identity element $0$. Here is a host of examples of isometric representations of $P$. Let $A \subset G$ be a non-empty set such that $P+A \subset A$.  Consider the Hilbert space $\ell^{2}(A)$ and let $\{\delta_{x}: x \in A\}$ be its standard orthonormal basis. For $a \in P$, let $V_a$ be the isometry on $\ell^2(A)$ defined by
\[
V_{a}(\delta_x):=\delta_{a+x}.\]
Then, $V^{A}:=\{V_{a}\}_{a \in P}$ is an isometric representation of $P$, which we call  the isometric representation associated to  $A$. Denote the $C^{*}$-algebra generated by $\{V_{a}: a \in P\}$ by $C^{*}(P,A)$. The case $A=P$ corresponds to the reduced $C^{*}$-algebra. The $C^{*}$-algebras $C^{*}(P,A)$ were analysed in great detail in \cite{Salas}  when $P=\mathbb{N}^{k}$.

A common feature that the isometric representations $V^{A}$ share, when we vary  $A$, is that the range projections $\{V_{a}V_{a}^{*}:a \in P\}$ form a commuting family. Thus, it is natural to consider  the following universal $C^{*}$-algebra. 
Let $C^{*}_{c}(P)$ be the universal  $C^{*}$-algebra generated by isometries $\{v_{a}:a \in P\}$ such that 
\begin{enumerate}
\item[(C1)] for $a,b \in P$, $v_{a}v_{b}=v_{a+b}$, and
\item[(C2)] for $a,b \in P$, $e_{a}e_{b}=e_{b}e_{a}$. Here, $e_{a}:=v_{a}v_{a}^{*}$ and $e_b:=v_bv_b^{*}$.
\end{enumerate}
The subscript 'c' in  $C_{c}^{*}(P)$ stands to  indicate that the range projections commute. 

The $C^{*}$-algebra $C_{c}^{*}(\bbn^k)$ was analysed from a groupoid perspective by Salas in \cite{Salas}. Murphy, in \cite{Murphy} by independent methods, studied $C_{c}^{*}(\bbn^2)$ and also the universal one generated by two commuting isometries. He proved that the latter $C^{*}$-algebra is complicated, in particular not nuclear, while the former is nuclear and more manageable.  The $C^{*}$-algebra $C_{c}^{*}(P)$ for a numerical semigroup $P$ was considered in \cite{Raeburn_Sean} and in \cite{Sean}. The analog of $C_{c}^{*}(P)$ in the topological setting  was investigated using groupoid methods by the last author in \cite{Sundar_Ore}. Both \cite{Salas} and \cite{Sundar_Ore} borrow substantial amount of material from \cite{Renault_Muhly}.

This paper is an attempt to understand the structure of KMS states on $C_{c}^{*}(P)$ for a natural time evolution. Let $c:G \to \mathbb{R}$ be a non-zero homomorphism. By the universal property of $C_{c}^{*}(P)$, for every $t \in \mathbb{R}$, there exists a $^*$-homomorphism $\sigma_{t}:C_{c}^{*}(P) \to C_{c}^{*}(P)$ such that 
\[
\sigma_{t}(v_a)=e^{itc(a)}v_{a}\]
for $a \in P$.  Then, $\sigma^{c}:=\{\sigma_{t}\}_{t \in \mathbb{R}}$ defines an action of $\mathbb{R}$ on $C_{c}^{*}(P)$. 

In this paper, we analyse the structure of KMS states for $\sigma^c$ for a toy model by discussing the case when $P=\mathbb{N}^{2}$.  We show in particular that,   in contrast to the reduced version $C_{red}^{*}(\mathbb{N}^{2})$, the set of KMS states on $C_{c}^{*}(\mathbb{N}^{2})$ has a  rich  structure.  We exhibit, in particular, uncountably many examples of extremal KMS states of type I, II and III. 
This paper is inspired by the discussions in \cite{Thomsen}. 

Next, we explain the results obtained.  Let $c:\mathbb{Z}^{2} \to \mathbb{R}$ be a non-zero homomorphism.  We normalise and assume that $c(1,0)=1$ and $c(0,1)=\theta$. Let $\Omega:=\{0,1\}^{\mathbb{Z}}$ be the Cantor space equipped with the product topology. Let $\tau:\Omega \to \Omega$ be the Bernoulli shift defined by
$
\tau(x)_k=x_{k-1}$.
Define $\chi: \Omega \to \mathbb{R}$ by \begin{equation*}
\label{isometries}
\chi(x):=\begin{cases}
 1  & \mbox{ if
} x_{-1} =0,\cr
   &\cr
   - \theta &  \mbox{ if } x_{-1}=1.
         \end{cases}
\end{equation*}
The results obtained in this paper are summarised below. 
\begin{enumerate}
\item[(1)] If either $\theta<0$ or $\beta<0$, there is no $\beta$-KMS state on $C_{c}^{*}(\mathbb{N}^{2})$. 
\item[(2)] Suppose $\theta=0$ and $\beta >0$. The simplex of $\beta$-KMS states is homeomorphic to the set of probability measures on $\mathbb{T}$. In this case, the extremal $\beta$-KMS states are of type I. 
\item[(3)] Suppose $\theta$ is irrational, positive and $\beta>0$. The simplex of $\beta$-KMS states is homeomorphic to the set of probability measures $m$ on the Cantor space $\Omega$ that are $e^{-\beta \chi}$-conformal, i.e. \[
m(\tau(B))=\int_B e^{-\beta \chi} dm\]
for every Borel subset $B$ of $\Omega$. 
\item[(4)] Suppose $\theta>0$ and $\beta>0$. For each $t \in \{I,II,III\}$, there  are uncountably many extremal $\beta$-KMS states of type $t$. Moreover, for $t \in \{II,III\}$, the extremal $\beta$-KMS states of type $t$ are in bijective correspondence with non-atomic type $t$ ergodic measures $m$ on $\Omega$ that are $e^{-\beta \chi}$-conformal.
\item[(5)] The simplex of tracial states of $C_{c}^{*}(\mathbb{N}^{2})$ is homeomorphic to the simplex of probability measures on the torus $\mathbb{T}^{2}$. 
\end{enumerate}

We end this introduction by mentioning a result obtained in Section 4, which we believe is worth highlighting.  Let $X$ be a compact metric space and let $\phi:X \to X$ be a homeomorphism. A probability measure $m$ on $X$ is said to be a conformal measure with potential $F$, where $F$ is a measurable function on $X$,
if 
\[
\frac{d(m \circ \phi)}{dm}=e^{F}.\]
Up to the authors' knowledge, the  reference in the literature for the proof of the existence of an ergodic type III conformal measure with continuous potential $F$  is the work of Katznelson (\cite{Katz}). (The reader is referred to Theorem  3.2, Theorem 3.3 of Part II in \cite{Katz} and is also referred to Theorem 5.2 of \cite{Thomsen}.  Also, as we will see in Section 5 that, Nakada's examples in \cite{Nakada} and in \cite{Nakada1} provide such examples).   If we demand only the measurability of the potential, then it is well known that odometers provide a rich source of such type III examples (see \cite{Dan_silva}). 

Here, by making use of Arnold's dyadic adding machine, we produce an example of an ergodic type III conformal measure with continuous potential on the Cantor space $\{0,1\}^{\mathbb{Z}}$, where the action is the usual shift. This construction is probably simpler, modulo accepting the fact that Arnold's adding machine is one of the simplest type III examples. But a drawback with our construction is that, unlike the example due to Katznelson, our dynamical system is not minimal. It is not clear to the authors whether this construction can be tweaked to produce an example, based on odometers, which is minimal.

\section{A groupoid model for $C_{c}^{*}(P)$}
In this section, we review the groupoid model for $C_{c}^{*}(P)$ described in \cite{Salas} and in \cite{Sundar_Ore}. We follow the exposition given in \cite{Sundar_Ore}. Let $G$ be a countable, discrete abelian group, and let $P$ be a 
subsemigroup of $G$ containing the identity element $0$. We assume that $P$ generates $G$, i.e. $P-P=G$. 

Recall that $C_{c}^{*}(P)$ is the universal unital $C^{*}$-algebra generated by a family of isometries $\{v_a: a \in P\}$ such that 
\begin{enumerate}
\item[(C1)] for $a,b \in P$, $v_av_b=v_{a+b}$, and
\item[(C2)] for $a,b \in P$, $e_ae_b=e_be_a$, where $e_a:=v_av_a^{*}$ and $e_b:=v_bv_b^{*}$. 
\end{enumerate}
Let $s \in G$ be given. Choose $a,b\in P$ such that $s=a-b$. Set $w_s:=v_b^{*}v_a$. Thanks to Prop. 3.4 of \cite{Sundar_Ore}, $w_s$ is well defined and  $\{w_s: s \in G\}$ is a family of partial isometries whose range projections
commute. For $s \in G$, let $e_s:=w_sw_s^*$. 

Next, we describe a groupoid whose $C^{*}$-algebra is a quotient of $C_{c}^{*}(P)$. Let $\mathcal{P}(G)$ be the power set of $G$ which we  identify, in the usual way, with $\{0,1\}^{G}$. Endow $\mathcal{P}(G)$ with the product topology inherited via this identification. Then, $\mathcal{P}(G)$ is a compact Hausdorff space and is metrisable. The map 
\[
\mathcal{P}(G) \times G \ni (A,s) \to A+s \in \mathcal{P}(G)\]
defines an action of $G$ on $\mathcal{P}(G)$. 

Let \[\overline{X}_{u}:=\{A \in \mathcal{P}(G): 0 \in A, -P+A \subset A\}.\] Note that $\overline{X}_u$ is a closed subset of $\mathcal{P}(G)$ and hence a compact subset of $\mathcal{P}(G)$. Also, $\overline{X}_u$ is $P$-invariant, i.e. if $A \in \overline{X}_u$ and $a \in P$, then $A+a \in \overline{X}_u$. Let $\mathcal{G}$ be the reduction of the transformation groupoid $\mathcal{P}(G) \rtimes G$ onto $\overline{X}_u$, i.e.
\[
\mathcal{G}:=\{(A,s) \in \mathcal{P}(G) \times G: A \in \overline{X}_u, A+s \in \overline{X}_u\}=\{(A,s): A \in \overline{X}_u, -s \in A\}.\] The multiplication and inversion on $\clg$ are given by
\begin{align*}
(A,s)(B,t)&:=(A,s+t) ~~~\textrm{~~~~~if $A+s=B$}, \textrm{~~and}\\
(A,s)^{-1}&:=(A+s,-s).
\end{align*}
The groupoid $\mathcal{G}$ is $r$-discrete and is the usual Deaconu-Renault groupoid $\overline{X}_u \rtimes P$.

For $f \in C(\overline{X}_u)$, define $\widetilde{f} \in C_{c}(\overline{X}_u \rtimes P)$ by \begin{equation}
\widetilde{f}(A,s):=\begin{cases}
 f(A)  & \mbox{ if
} s =0, \cr
   &\cr
    0 &  \mbox{ if $s \neq 0$. }
         \end{cases}
\end{equation}
Observe that $C(\overline{X}_u) \ni f \to \widetilde{f} \in C_{c}(\overline{X}_u \rtimes P)$ is a unital $^*$-algebra homomorphism which is injective. Via this embedding, we identify $C(\overline{X}_u)$ as a unital $^*$-subalgebra of $C_{c}(\mathcal{G})$.  For $f \in C(\overline{X}_u)$, we abuse notation and we denote $\widetilde{f}$ by $f$. 

For $s \in G$, let $\overline{w}_s \in C_{c}(\overline{X}_u \rtimes P)$ be defined by \begin{equation}
\overline{w}_s(A,t):=\begin{cases}
 1  & \mbox{ if
} t=-s, \cr
   &\cr
    0 &  \mbox{ if $t \neq -s$}.
         \end{cases}
\end{equation}
For $a \in P$, set $\overline{v}_a:=\overline{w}_a$.   Observe the following facts.  The details involving routine computations are left to the reader. 
\begin{enumerate}
\item[(1)] For $a \in P$, $\overline{v}_a$ is an isometry and $\overline{v}_a\overline{v}_b=\overline{v}_{a+b}$. 
\item[(2)] For $a,b \in P$ and $s=a-b$, $\overline{w}_s=\overline{v}_b^{*}\overline{v}_a$. 
\item[(3)] For $s \in G$, let $\epsilon_{s} \in C(\overline{X}_u)$ be defined by the equation $\epsilon_{s}(A)=1_{A}(s)$. Then, $\overline{w}_s\overline{w}_s^{*}=\epsilon_{s}$ for every $s \in G$. Consequently,  $\{\overline{w}_s\overline{w}_s^{*}:s \in G\}$ generates $C(\overline{X}_u)$. Also, the range projections $\{\overline{v}_a\overline{v}_a^{*}:a \in P\}$ form a commuting family. 
\item[(4)] The $C^{*}$-algebra generated by $\{\overline{v}_a:a \in P\}$ is $C^{*}(\mathcal{G})$. To see this, consider a function $f \in C_{c}(\clg)=C_c(\overline{X}_u \rtimes P)$. It suffices to consider the case when $f$ is supported on $\overline{X}_u \times\{-s\}$ for some $s \in G$. Let $h:\overline{X}_u \to \mathbb{C}$ be defined by 
\begin{equation*}
h(A):=\begin{cases}
 f(A,-s)  & \mbox{ if
} (A,-s) \in \overline{X}_u \rtimes P, \cr
   &\cr
    0 &  \mbox{ otherwise}.
         \end{cases}
\end{equation*}
Then, $h$ is continuous and $h*\overline{w}_s=f$. By $(3)$, $h$ lies in the $C^{*}$-algebra generated by $\{\overline{w}_s:s \in G\}$. Therefore, $\{\overline{w}_s:s \in G\}$ generates $C^{*}(\mathcal{G})$. It follows from $(2)$ that $\{\overline{v}_a: a \in P\}$ generates $C^{*}(\overline{X}_u \rtimes P)$. 
\end{enumerate}
Hence, there exists a unique surjective $^*$-homomorphism $\Phi:C_{c}^{*}(P) \to C^{*}(\overline{X}_u \rtimes P)$ such that \begin{equation}
\label{main equation}
\Phi(v_a)=\overline{v}_a.\end{equation}
Theorem 7.4 of \cite{Sundar_Ore} asserts that $\Phi$ is invertible and $\Phi$ is an isomorphism. 
Via the map $\Phi$, we identify $C_{c}^{*}(P)$ with $C^{*}(\overline{X}_u \rtimes P)$.  Henceforth, we do not distinguish between $w_s$ and $\overline{w}_s$. 

For $f \in C(\overline{X}_u)$ and $s \in G$, let $R_{s}(f) \in C(\overline{X}_u)$ be defined by the equation
\begin{equation*}
R_s(f)(A):=\begin{cases}
 f(A-s)  & \mbox{ if~} A-s \in \overline{X}_u, \cr
   &\cr
    0 &  \mbox{ otherwise}.
         \end{cases}
\end{equation*}
For $s \in G$ and $f \in C(\overline{X}_u) \subset C^{*}(\overline{X}_u \rtimes P)$, the  covariance relation
\[
w_s f w_s^{*}=R_s(f)\]
is satisfied.
A pleasant consequence of the above covariance relation is the fact that  $\textrm{span}\{fw_s: f \in C(\overline{X}_u), s \in G\}$ is a unital dense $^*$-subalgebra of $C^{*}(\overline{X}_u \rtimes P)$.

Let \[\overline{Y}_{u}:=\{A \in \mathcal{P}(G): A \neq \emptyset, -P+A \subset A\}.\]
Observe that $\overline{Y}_u$ is a locally compact Hausdorff space and $G$ leaves $\overline{Y}_u$ invariant. Also, $\overline{X}_u$ is a clopen set in $\overline{Y}_u$.  It is clear that the groupoid $\overline{X}_u \rtimes P$ is the reduction of the 
transformation groupoid $\overline{Y}_u \rtimes G$ onto $\overline{X}_u$. This has the consequence that $C^{*}(\overline{X}_u \rtimes P)$ is isomorphic to the cut-down $p(C_0(\overline{Y}_u) \rtimes G)p$ where $p=1_{\overline{X}_u}$. 
  Moreover,  the union $\displaystyle \bigcup_{a \in P}(\overline{X}_u-a)=\overline{Y}_u$. Consequently, it follows that $p$ is a full projection in $C_0(\overline{Y}_u) \rtimes G$. 

  Often, we embedd $C_{c}^{*}(P)$ inside the crossed product $C_0(\overline{Y}_u) \rtimes G$ and view $C_{c}^{*}(P)$ as a full corner of $C_{0}(\overline{Y}_u) \rtimes G$. The embedding $C_{c}^{*}(P) \to C_0(\overline{Y}_u) \rtimes G$ is 
  given by the rules \[C(\overline{X}_u) \ni f \to f1_{\overline{X}_u} \in C_0(\overline{Y}_u) \textrm{~and~} w_s \to u_s1_{\overline{X}_u}.\] 
  Here, $\{u_s: s \in G\}$ are the canonical unitaries of the crossed product $C_0(\overline{Y}_u)\rtimes G$.
  
 \textbf{Tracial states on $C_{c}^{*}(P)$:}  We first determine the tracial states on $C_c^{*}(P)$. 
 Consider the group $C^{*}$-algebra $C^{*}(G)$, and let $\{u_s:s \in G\}$ be the canonical unitaries of $C^{*}(G)$. By the universal property, there exists a surjective $^*$-homomorphism 
 $\pi:C_c^{*}(P) \to C^{*}(G)$ such that $\pi(v_a)=u_a$. In the groupoid picture, $\pi$ coincides with the restriction map $Res:C^{*}(\overline{X}_u \rtimes P) \to C^{*}(\overline{X}_u \rtimes P|_{\{G\}}) \cong C^{*}(G)$. Note that $\{G\}$ is an invariant
 closed subset of the unit space $\overline{X}_u$ of the groupoid $\mathcal{G}=\overline{X}_u \rtimes P$.

  \begin{ppsn}
  \label{tracial}
  Every tracial state on $C_{c}^{*}(P)$ factors through $C^{*}(G)$. Thus, tracial states on $C_{c}^{*}(P)$ are in bijective correspondence with probability measures on $\widehat{G}$. 
    \end{ppsn}
  \textit{Proof.} Let $F:=\{G\}$ and $X_u:=\overline{X}_u\backslash F$. 
    Let $\omega$ be a tracial state on $C_{c}^{*}(P)$. Let $m$ be the measure on $\overline{X}_u$ that corresponds to the state $\omega|_{C(\overline{X}_u)}$.
  Let $a \in P$ be given. Note that $v_av_a^{*}=1_{\overline{X}_u+a}$. Since, $\omega$ is tracial, we have
  \[
  m(\overline{X}_u+a)=\omega(v_av_a^*)=\omega(v_a^*v_a)=\omega(1)=1=m(\overline{X}_u).\]
  But, $\overline{X}_u+a \subset \overline{X}_u$. Hence, $\overline{X}_u \backslash (\overline{X}_u+a)$ has measure zero. Consequently, the set $X_u=\bigcup_{a \in P}\overline{X}_u \backslash (\overline{X}_u+a)$ has measure zero. 
  Therefore, $m$ is concentrated on $F=\{G\}$. 
  
  Since $X_u$ is an open invariant set, we have the following short-exact sequence
  \[
0 \longrightarrow C^{*}(\mathcal{G}|_{X_u}) \longrightarrow C^{*}(\mathcal{G}) \to C^{*}(\mathcal{G}|_F) \cong C^{*}(G) \longrightarrow 0.\]

To complete the proof, it suffices to show that $\omega$ vanishes on $C^{*}(\mathcal{G}|_{X_u})$. 
However, $span\{fw_s: f \in C_c(X_u), s \in G\}$ is a dense $^*$-subalgebra of $C^{*}(\mathcal{G}|_{X_u})$. Thus, it is enough to prove that $\omega(fw_s)=0$ whenever $f \in C_c(X_u)$ and $s \in G$. 
Let $f \in C_c(X_u)$ and $s \in G$ be given. The equality $\omega(fw_s)=0$ follows from the Cauchy-Schwartz inequality and the fact that $m(X_u)=0$. This completes the proof. \hfill $\Box$
  
  We end this section with  a few definitions that we need later.
  \begin{dfn}
 Let $(Y,\mathcal{B})$ be a standard Borel space on which $G$ acts measurably.  Let $X$ be a Borel subset of $Y$. We say that $(Y,X)$ is a $(G,P)$-space if 
\begin{enumerate}
\item[(1)] the set $X$ is  $P$-invariant, i.e. $X+P \subset X$,
 and
\item[(2)] the union $\displaystyle \bigcup_{a \in P}(X-a)=Y$. 
\end{enumerate}
Let $(Y,X)$ be a $(G,P)$-space. We say that $(Y,X)$ is pure if  $\displaystyle \bigcap_{a \in P}(X+a)=\emptyset$.
    \end{dfn}
    
    \begin{xmpl}
    The pair $(\overline{Y}_u,\overline{X}_u)$ is a $(G,P)$-space. Define \[Y_u:=\overline{Y}_u \backslash \{G\};~~X_u:= \overline{X}_u \backslash \{G\}.\] 
    Then, $(Y_u,X_u)$ is pure. 
    
    \end{xmpl}
    
    \begin{dfn}
Let $(Y,X)$ be a $(G,P)$-space. Let $c:G \to \mathbb{R}$ be a non-zero homomorphism and let $\beta$ be a fixed real number. 
Suppose $m$ is a measure on $Y$. We say that $m$ is an $e^{-\beta c}$-conformal measure on the $(G,P)$-space $(Y,X)$ if
\begin{enumerate}
\item[(i)]  $m(X)=1$, and
\item[(ii)] for every Borel set $E \subset Y$ and 
$s \in G$, 
\[
m(E+s)=e^{-\beta c(s)}m(E).\]
\end{enumerate}
\end{dfn}
We often abuse terminology and call an $e^{-\beta c}$-conformal measure on a $(G,P)$-space $(Y,X)$ an $e^{-\beta c}$-conformal measure on $Y$.

  \section{KMS states and conformal measures}
  
  Let us recall the definition of a  $\beta$-KMS state.  Let $A$ be a $C^*$-algebra and suppose $\tau:=\{\tau_t\}_{t \in \bbr}$ is a $1$-parameter group of automorphisms of $A$.  Suppose $\beta \in \mathbb{R}$. A state $\omega$ on $A$ is a $\beta$-KMS state for $\tau$ if 
	$$\omega(ab)=\omega(b\tau_{i\beta}(a))$$
	for all $a, b$ in a norm dense $\tau$-invariant $^*$-algebra of analytic elements in $A$. 
	
	Let $\omega$ be a $\beta$-KMS state on $A$. Then, $\omega$ is said to be \emph{extremal} if it is an extreme point in the simplex of $\beta$-KMS states. It is well known that a $\beta$-KMS state is 
	extremal if and only if the associated GNS representation $\pi_\omega$ is factorial, i.e. $\pi_\omega(A)^{''}$ is a factor. We say an extremal $\beta$-KMS state is of type $t$ if $\pi_\omega(A)^{''}$ is a factor
	of type $t$. 
	
	Fix a non-zero homomorphism $c:G \to \bbr$. Recall the $1$-parameter group of automorphisms $\sigma^{c}:=\{\sigma_{t}\}_{t \in \mathbb{R}}$ defined on $C_{c}^{*}(P)\cong C^{*}(\overline{X}_u \rtimes P)$ by  
	\[
	\sigma_{t}(v_a)=e^{itc(a)}v_{a}\]
	for $a \in P$.
For $t \in G$, let $D(t), R(t)\subset \overline{X}_u$ be defined by \begin{align*}
D(t)&:=\{A\in \overline{X}_u: A+t \in \overline{X}_u \} \\
R(t)&:=\{A\in \overline{X}_u: A-t \in \overline{X}_u \}.
\end{align*}
For $t \in G$, denote the map $D(t) \ni A\to A+t \in R(t)$ by $T_t$. 

Let $\beta$ be a real number. As we have already described the tracial states on $C_c^{*}(P)$  in Prop. \ref{tracial}, we assume for the rest of this section that $\beta$ is an arbitrary, but a fixed non-zero real number.  Let $\omega$ be a $\beta$-KMS state on the $C^{*}$-algebra $C_c^{*}(P)=C^*(\overline{X}_u \rtimes P)$ for $\sigma^c$. Restriction of $\omega$ to $C(\overline{X}_u)$ defines a  probability measure $m$ on $\overline{X}_u$. It can be checked  from the covariance relation and the KMS condition that \[\omega(f)=e^{-\beta c(t)}\omega(f\circ T_t)\] for every $f \in C(R(t))$, or equivalently $m(E+t)=e^{-\beta c(t)}m(E)$ for every Borel subset $E$ of $ D(t) $. 
 
 For $s,t \in G$, we say $s \leq t$ if $t-s \in P$. Then,  $G$ with the preorder $\leq$ is a directed set. Note that $P$ is a cofinal subset of $G$. Moreover, if $s \leq t$, then  $\overline{X}_u-s\subseteq \overline{X}_u-t$. Let $s \in G$. Define a measure $m_s$ on  $\overline{X}_u-s$ by \[m_s(E)=e^{\beta c(s)}m(E+s).\] Thanks to the fact that $m$ is conformal, it is clear that $m_t|_{\overline{X}_u-s}=m_s$ for $s\leq t$. Define a $\sigma$-finite measure $m_\omega$ on $\overline{Y}_u$ by setting
\begin{equation}\label{conformalmeasure}
m_\omega(E):= \lim_{s \in G}m_s(E \cap (\overline{X}_u-s))=\lim_{s \in P}m_s(E\cap (\overline{X}_u-s)). 
\end{equation}
\begin{ppsn}{\label{statetomeasure}}
	The measure $m_\omega$ on $\overline{Y}_u$ defined as in Equation \ref{conformalmeasure} is an $e^{-\beta c}$-conformal measure on the $(G,P)$-space $(\overline{Y}_u,\overline{X}_u)$.
\end{ppsn}
\textit{Proof.}
It is clear that $m_\omega|_{\overline{X}_u}=m$ and so $m_\omega(\overline{X}_u)=1$.
 For $t\in G$ and  a Borel subset $E\subset \overline{Y}_u$, 
\begin{eqnarray*}
m_\omega(E+t)&=&\lim_sm_s((E+t)\cap (\overline{X}_u-s))\\
&=&\lim_s e^{\beta c(s)}m(((E+t)\cap (\overline{X}_u-s))+s)\\
&=& \lim_s e^{\beta c(s)}m((E\cap( \overline{X}_u-(s+t)))+s+t)\\
&=& e^{-\beta c(t)}\lim_r e^{\beta c(r)} m((E\cap (\overline{X}_u-r))+r)~~(\textrm{by a change of variable $r=s+t$})\\
&=& e^{-\beta c(t)}\lim_r m_r(E\cap (\overline{X}_u-r))\\
&=& e^{-\beta c(t)}m_\omega(E).
\end{eqnarray*}
Hence the proof. \hfill $\Box$

Recall that for an $r$-discrete groupoid $\clg$, the map $f \mapsto f|_{\clg^{(0)}}$ from $C_c(\clg)$ to $C_0(\clg^{(0)})$ defines a conditional expectation  $E:C^*(\clg) \to C_0(\clg^{(0)})$. For $\mathcal{G}:=\overline{X}_u \rtimes P$, the conditional expectation $E$ has the following form %
 \begin{equation*}
 E(fw_s):=\begin{cases}
 f  & \mbox{ if~~} s=0,\\
 0 &  \mbox{ otherwise}.
 \end{cases}
 \end{equation*}
 
\begin{ppsn}\label{measuretostate}
	The map $\omega \mapsto m_{\omega}$ from the set of $\beta$-KMS states for $\sigma=\sigma^c$ to the set of $e^{-\beta c}$-conformal measures on $\overline{Y}_u$ is surjective. More specifically, for an $e^{-\beta c}$-conformal $m$ on $\overline{Y}_u$, the state $\omega_m$ on $C^{*}(\overline{X}_u \rtimes P)$ defined by 
	\begin{equation*}\label{stateexp}
	\omega_m(a) = \int_{\overline{X}_u} E(a)~dm 
	\end{equation*}
is a $\beta$-KMS state such that $m_{\omega_{m}}=m$.
\end{ppsn}
\textit{Proof.}
Let $f,g \in C(\overline{X}_u)$ and let $s,t \in G$ be given. Suppose $s+t \neq 0$. It follows  from the covariance relation that $fw_sgw_t=fR_s(g)1_{\overline{X}_u+s}w_{s+t}$. Hence, $\omega_m(fw_sgw_t)=0$. By the same reasoning, we have 
$\omega_m(gw_t\sigma_{i\beta}(fw_s))=e^{-\beta c(s)}\omega_m(gw_t fw_s)=0$. Thus, if $s+t \neq 0$, then 
$
\omega_m(fw_sgw_t)=\omega_m(gw_t \sigma_{i\beta}(fw_s))$.

Assume now that $s+t=0$, i.e $t=-s$. Then, 
\begin{equation}\label{kmsequation}
\omega_m(fw_sgw_{-s})=w_m(fR_s(g))=\int_{R(s)}fR_s(g)~dm
\end{equation}
since $R_s(g)$ vanishes outside $R(s)$. On the other hand,
\begin{eqnarray*}
	\omega_m(gw_{-s}\sigma_{i\beta}(fw_s))&=& e^{-\beta c(s)}\omega_m(gw_{-s}fw_s)\\
	&=& e^{-\beta c(s)}\int_{D(s)}R_{-s}(f)g ~dm \qquad (\text{since}~R_{-s}(f) \text{~vanishes outside}~ D(s))\\
	&= &  e^{-\beta c(s)}\int_{D(s)}(fR_s(g))\circ T_s ~dm\\
	&=& e^{-\beta c(s)}\int_{R(s)}fR_s(g)d(T_s)_*m\\
&=& e^{-\beta c(s)}e^{\beta c(s)}\int_{R(s)}fR_s(g)~dm\qquad (\text{since}~\frac{d(T_s)_*m}{dm}=e^{\beta c(s)}~~ \text{on}~R(s))\\
&=& \omega_m(fw_sgw_{-s})\qquad\qquad \qquad\qquad (\text{by Eq.}~\ref{kmsequation}).
\end{eqnarray*}
Hence, $\omega_m$ is a $\beta$-KMS state. Clearly $m_{\omega_{m}}=m$. The proof is complete
\hfill $\Box$

\begin{rmrk}
In general, the map  $\omega \to m_\omega$ of Prop. \ref{statetomeasure} need not be injective. The structure of a generic KMS state is determined by Neshveyev's result (Thm. 1.3, \cite{Neshveyev}). 
Let $\mathcal{G}$ be an $r$-discrete groupoid and let $\overline{c}:\mathcal{G} \to \mathbb{R}$ be a continuous homomorphism. For $t \in \bbr$, let $\sigma_t$ be the automorphism of $C^{*}(\clg)$ defined by
\[
\sigma_t(f)(\gamma):=e^{it\overline{c}(\gamma)}f(\gamma)\]
for $f \in C_{c}(\clg)$. Then, $\sigma^{\overline{c}}:=\{\sigma_t\}_{t \in \bbr}$ defines a $1$-parameter group of automorphisms on $C^{*}(\clg)$. 
Theorem 1.3 of \cite{Neshveyev} describes the set of $\beta$-KMS states for $\sigma^{\overline{c}}$. 

In our situation, the groupoid $\mathcal{G}$ is the Deaconu-Renault groupoid $\overline{X}_u \rtimes P$ and the homomorphism $\overline{c}:\mathcal{G} \to \bbr$ is given by 
\[
\overline{c}(A,g)=-c(g).\]
One immediate corollary of Theorem 1.3 of \cite{Neshveyev} is that, if  the homomorphism $c:G \to \bbr$ is injective, then the map $\omega \to m_{\omega}$ of Prop. \ref{statetomeasure}
is a bijection. 
\end{rmrk}

\begin{ppsn}
\label{ergodicity}
	Let $\beta \neq 0$ and let $\omega$ be an extremal $\beta$-KMS state for $\sigma=\sigma^{c}$ on the $C^{*}$-algebra $C^*(\overline{X}_u \rtimes P)$. Then, the $e^{-\beta c}$-conformal measure $m_\omega$ on $\overline{Y}_u$ defined as in Proposition \ref{statetomeasure} is ergodic for the $G$-action on $\overline{Y}_u$.
\end{ppsn}
\textit{Proof.}
The proof is an adaptation of the proof of Lemma 3.6 in \cite{Thomsen} with appropriate modifications. 
Assume that $m_\omega$ is not ergodic. Then, there exists a $G$-invariant Borel subset $A\subset \overline{Y}_u$ such that $m_\omega(A) \neq 0$ and $m_\omega(A^c) \neq 0$.

We claim that  $0<m_\omega(A\cap \overline{X}_u)<1$ and 
$0<m_\omega(A^c\cap  \overline{X}_u)<1$.

Suppose that $m_\omega(A\cap \overline{X}_u)=0$.   By conformality,  $m_\omega((A\cap \overline{X}_u)-t)=0$ for every $t \in G$. Hence, $\bigcup_{t \in G}(A\cap \overline{X}_u)-t$ is a null set. The fact that $A$ is $G$-invariant implies that  $\bigcup_{t \in G}(A\cap \overline{X}_u)-t=A \cap \big(\bigcup_{t \in G}(\overline{X}_u-t)\big)=A$. Therefore, $m_\omega(A)=0$ which is a contradiction. This proves that 
$m_\omega(A\cap \overline{X}_u)\neq 0$. Similarly, $m_\omega(A^c \cap \overline{X}_u) \neq 0$. Consequently, $0 <m_\omega(A \cap \overline{X}_u)<1$ and $0 < m_\omega(A^c \cap \overline{X}_u)<1$. This proves the claim.

 Let $(H_\omega,\pi_\omega, \Omega_\omega) $ be the GNS representation of $\omega$. Extend $\pi_\omega|_{C(\overline{X}_u)}$  to a $^*$- homomorphism $\overline{\pi}_\omega: L^{\infty}(\overline{X}_u,m_\omega) \to \pi_\omega(C^*(\overline{X}_u \rtimes P))''$ by defining \[\overline{\pi}_\omega(h):=\lim_{ n \to \infty}\pi_\omega(f_n)\] 
 Here, the limit is taken in the SOT sense and $\{f_n\}$ is any sequence in $C(\overline{X}_u)$ such that for every $n$,  $|f_n| \leq ||h||_\infty$ a.e and the sequence $(f_n(x)) \to h(x)$ for almost all $x$. 
  
  For $s \in G$, set  $W_s=\pi_\omega(w_s)$. Then, the following covariance relation is satisfied 
\[ W_s\overline{\pi}_\omega(f) W_s^*=\overline{\pi}_\omega(R_s(f))\]
for $s \in G$ and $f \in L^{\infty}(\overline{X}_u,m_\omega)$.

Next, we claim that $\overline{\pi}_\omega(1_{A \cap \overline{X}_u})$ and $\overline{\pi}_\omega(1_{A^c \cap \overline{X}_u})$ are central in $\pi_\omega(C^*(\overline{X}_u \rtimes P))''$. 
We will give details for $\overline{\pi}_\omega(1_{A \cap \overline{X}_u})$ and the centrality of $\overline{\pi}_\omega(1_{A^c \cap \overline{X}_u})$ follows similarly. 
Let $s \in P$ be given. Note that $R_{-s}(1_{A \cap \overline{X}_u})=1_{(A-s) \cap (\overline{X}_u-s)\cap \overline{X}_u}=1_{A \cap \overline{X}_u}$ and $\overline{\pi}_\omega(f)$ commutes with $W_sW_s^*$ for every $f \in L^{\infty}(\overline{X}_u,m_\omega)$. 

Calculate as follows to observe that 
\begin{eqnarray*}
\overline{\pi}_\omega(1_{A \cap \overline{X}_u})W_s&=& \overline{\pi}_\omega(1_{A \cap \overline{X}_u})W_sW_s^*W_s\\
&=& W_sW_s^*\overline{\pi}_\omega(1_{A \cap \overline{X}_u})W_s\\
&=&W_s \overline{\pi}_\omega(R_{-s}(1_{A \cap \overline{X}_u}))\\
&=&W_s \overline{\pi}_\omega(1_{A \cap \overline{X}_u}).
\end{eqnarray*}
Since $\{W_s: s \in P\}$ generates $(\pi_\omega(C^*(\overline{X}_u \rtimes P))^{''}$, we can conclude that $\pi_\omega(1_A \cap \overline{X}_u)$ is central. This proves the claim. 

Define states $\omega_1$ and $\omega_2$ on $C^{*}(\overline{X}_u \rtimes P)$ by   \begin{align*}
\omega_1(a)&=\frac{1}{m_\omega(A\cap \overline{X}_u)}\langle \overline{\pi}_\omega(1_{A \cap \overline{X}_u})\pi_\omega(a)\Omega_\omega,\Omega_\omega\rangle, \\
\omega_2(a)&=\frac{1}{m_\omega(A^c\cap \overline{X}_u)}\langle \overline{\pi}_\omega(1_{A^c \cap \overline{X}_u})\pi_\omega(a)\Omega_\omega,\Omega_\omega\rangle.
\end{align*}
The centrality of the projections $\overline{\pi_\omega}(1_{A \cap \overline{X}_u})$ and $\overline{\pi_\omega}(1_{A^c \cap \overline{X}_u})$ imply that $\omega_1$ and $\omega_2$ 
are $\beta$-KMS states on $C^*(\overline{X}_u \rtimes P)$. Moreover, $\omega=m_\omega(A\cap \overline{X}_u)\omega_1+m_\omega(A^c\cap \overline{X}_u)\omega_2$ which contradicts the extremality of $\omega$.
Hence the proof. 
\hfill $\Box$

Next, we work out the GNS representation of the KMS state $\omega_m$ given by Prop. \ref{measuretostate}. Suppose $m$ is an $e^{-\beta c}$ conformal measure on $\overline{Y}_u$ and let $\omega_m$ be the $\beta$-KMS state on the $C^*$-algebra $C_{c}^{*}(P)=C^*(\overline{X}_u \rtimes P)$ obtained via the conditional expectation as in Prop. \ref{measuretostate}. Our goal is to show that if $(H_\omega,\pi_\omega, \Omega_\omega)$ is the GNS  representation of $\omega_m$, then the von Neumann algebra  $(\pi_\omega(C^*(\overline{X}_u \rtimes P)))''$  is isomorphic to the full corner $1_{\overline{X}_u}
(L^\infty(\overline{Y}_u)\rtimes G)1_{\overline{X}_u}$ of $L^\infty(\overline{Y}_u)\rtimes G$ . 

Let $\lambda:=\{\lambda_s\}_{s \in G}$ be the Koopman representation of $G$ on $L^{2}(\overline{Y}_u,m)$. Recall that 
\[
\lambda_s\xi(A)=e^{\frac{\beta c(s)}{2}}\xi(A-s)\]
for $s \in G$ and $\xi \in L^{2}(\overline{Y}_u,m)$. 

Let $\clk=l^2(G,L^2(\overline{Y}_u))$ and $\pi_0$ be the  representation of $C_0(\overline{Y}_u)\rtimes G$ on $\clk$ defined by 
\begin{eqnarray*}
\pi_0(fu_t)(\xi)(s)=f\lambda_t(\xi(s-t)).
\end{eqnarray*}
Then, $(\pi_0(C_0(\overline{Y}_u) \rtimes G))''=L^\infty(\overline{Y}_u)\rtimes G$.  Let $\clh$ be the Hilbert subspace of $\clk$ given by \[\clh:=\{\xi \in l^2(G,L^2(\overline{Y}_u)): \xi(t)\in L^2(\overline{X}_u+t)\}.\] It is clear that $L^\infty(\overline{Y}_u)\rtimes G$ leaves $\clh$ invariant; thus giving a normal representation $\overline{\pi}$ of the crossed product $L^\infty(\overline{Y}_u)\rtimes G$ on $\clh$. Let $\pi=\overline{\pi}\circ\pi_0$. 

\textit{Notation:} For $s \in G$ and $\xi \in L^{2}(\overline{X}_u+s)$, let $\xi \otimes \delta_s \in \clh$ be defined by   \begin{equation*}
 \xi \otimes \delta_s(t)=\begin{cases}
 \xi  & \mbox{ if~~} t=s,\\
 0 &  \mbox{ otherwise}.
 \end{cases}
 \end{equation*}

\begin{lmma}\label{vonneumannalgebra}
	The von Neumann algebra  $(\pi(C_0(\overline{Y}_u) \rtimes G))'' $ is isomorphic to $L^\infty(\overline{Y}_u)\rtimes G$.
\end{lmma}
\textit{Proof.}
Note that since $\pi=\overline{\pi}\circ\pi_0$ and the von Neumann algebra generated by $\pi_0(C_0(\overline{Y}_u) \rtimes G)$ is $L^\infty(\overline{Y}_u)\rtimes G$, to complete the proof, it is enough to prove that $\overline{\pi}$ is injective.

First, we claim that $\overline{\pi}$ restricted to $L^{\infty}(\overline{Y}_u)$ is injective. Let $f \in L^{\infty}(\overline{Y}_u)$ be such that $\overline{\pi}(f)=0$. Suppose $s \in G$. Then, for every $\xi \in L^{2}(\overline{X}_u+s)$, 
\[
0=\overline{\pi}(f)(\xi \otimes \delta_s)(s)=f\xi.\]
Thus, $f=0$ a.e on $\overline{X}_u+s$ for every $s\in G$.  Since $\overline{Y}_u=\bigcup_{s \in G} \overline{X}_u+s$, it follows that $f=0$ a.e on $\overline{Y}_u$ proving that $\overline{\pi}$ is faithful on $L^\infty(\overline{Y}_u)$.

Let $E$ be the usual conditional expectation from $L^\infty(\overline{Y}_u)\rtimes G\to L^\infty(\overline{Y}_u)$ given by 
$
E(fu_s)=\delta_{s,0}f$. Here, $\delta_{s,0}$ is the Kronecker delta. Then, the conditional expectation $E$ is faithful. 
Write $\displaystyle \clh=\bigoplus_{t \in G}L^2(\overline{X}_u+t)$ and for $ t \in G$, let $P_t$ be the projection onto the subspace $L^2(\overline{X}_u+t)$. Define $\widetilde{E}: B(\clh) \to B(\clh)$ by \[\widetilde{E}(T)=\sum_{t \in G}P_tTP_t\] where the convergent sum is w.r.t.  the strong operator topology. Clearly, $\widetilde{E}\circ \overline{\pi}=\overline{\pi}\circ E$. 

Let $x \in L^{\infty}(\overline{Y}_u) \rtimes G$ be such that $\overline{\pi}(x)=0$. Then, $\overline{\pi}(x^*x)=0$. Therefore, 
\[
\overline{\pi}(E(x^*x))=\widetilde{E}(\overline{\pi}(x^*x))=0.\]
Since $\overline{\pi}$ is faithful on $L^{\infty}(\overline{Y}_u)$, we have $E(x^*x)=0$. Since $E$ is faithful, $x=0$. This completes the proof. \hfill $\Box$

 Recall that $C^{*}(\overline{X}_u \rtimes P)$ is  the full corner $p(C_0(\overline{Y}_u) \rtimes G)p$ where $p=1_{\overline{X}_u}$. 
 Consider the Hilbert space \[P\clh=\{\xi \in l^2(G,L^2(\overline{Y}_u)): \xi(t)\in L^2(\overline{X}_u\cap( \overline{X}_u+t))\}\] where $P=\pi(p)$. Define $\widetilde{\pi}: C^{*}(\overline{X}_u \rtimes P) \to B(P\clh)$ by $\widetilde{\pi}(fw_t)=P\pi(\tilde{f}u_t)P$ where $\tilde{f}=f1_{\overline{X}_u}$. Define $\xi \in P\clh$  by,
  \begin{equation*}
 \xi(t)=\begin{cases}
 1_{\overline{X}_u}  & \mbox{ if~~} t=0,\\
 0 &  \mbox{ otherwise}.
 \end{cases}
 \end{equation*}

\begin{lmma}{\label{gnsrepresentation}}
	With the foregoing notation, the triple $(\widetilde{\pi},P\clh,\xi)$ is  the GNS-representation of the $\beta$-KMS state $\omega_m$.
	Moreover, the von Neumann algebra $\widetilde{\pi}(C^{*}(\overline{X}_u \rtimes P))^{''}$ is isomorphic to the full corner $1_{\overline{X}_u}(L^{\infty}(\overline{Y}_u) \rtimes G)1_{\overline{X}_u}$. 
	\end{lmma}
\textit{Proof.} 
 Observe that, for $f \in C(\overline{X_u})$,  and $s,t\in G$,
 \begin{eqnarray*}
 	\widetilde{\pi}(fw_t)\xi(s)&=&P\pi(\tilde f u_t)P\xi(s)\\
 	&=&\pi(1_{\overline{X}_u}\tilde f u_t )\xi(s)\\
 	&=& 1_{\overline{X}_u}\tilde f \lambda_t(\xi(s-t))\\
 	&=& \delta_{s,t}1_{\overline{X}_u}f \lambda_t(\xi(0))\\
 	&=&e^{\frac{\beta c(s)}{2}}\delta_{s,t} 1_{\overline{X}_u}f1_{\overline{X}_u+s}
 \end{eqnarray*}
 From the above formula,  it is clear that $\langle \widetilde{\pi}(a)\xi,\xi\rangle =\omega_m(a)$ for every $a \in C^{*}(\overline{X}_u \rtimes P)$. 
Also, $\{\widetilde{\pi}(fw_t)\xi: f \in C(\overline{X}_u), t \in G\}=span\{\xi \otimes \delta_s: \xi \in C(\overline{X}_u \cap (\overline{X}_u+s)), s \in G\}$ and the latter set is dense in $P\clh$ proving that $\xi$ is cyclic for $\widetilde{\pi}$. 
Thus, $(\widetilde{\pi},P\clh,\xi)$ is the GNS representation of the $\beta$-KMS state $\omega_m$. 

From the definition of $\widetilde{\pi}$, we have
\[(\widetilde{\pi}(C^*(\overline{X}_u \rtimes P)))''=(P\pi(C_0(\overline{Y}_u) \rtimes G)P)''=P(\pi(C_0(\overline{Y}_u) \rtimes G))''P.\] The conclusion is now clear from Lemma \ref{vonneumannalgebra}.



We have the following main theorem of this section establishing the factor types of extremal $\beta$-KMS states $\omega_m$. 
\begin{thm}
\label{main_sec3}
	Suppose $m$ is an $e^{-\beta c}$-conformal measure on $\overline{Y}_u$. Let $\omega:=\omega_m$ be the $\beta$-KMS state as in Prop. \ref{measuretostate}. Denote   the GNS representation of $\omega$ by $\pi_\omega$. Suppose $t \in \{I_\infty,II_\infty,III\}$.  Then,
	\begin{itemize}
		\item[(1)]
		$(\pi_\omega(C^*(\overline{X}_u \rtimes P)))''$ is isomorphic to to the full corner $1_{\overline{X}_u}
		(L^\infty(\overline{Y}_u)\rtimes G)1_{\overline{X}_u} $.
		\item[(2)] $\omega_m$ is extremal $\iff L^\infty(\overline{Y}_u)\rtimes G$ is a factor $\iff m$ is ergodic and the $G$-action is free. Moreover,
		$\omega_m$ is of type $t$ if and only if $m$ is of type $t$.
		 	\end{itemize}
\end{thm}
\textit{Proof.} The first equivalence in $(2)$ follows from Lemma \ref{gnsrepresentation} and the  following facts. It is well known that 
\begin{enumerate}
\item[(i)] a $\beta$-KMS state $\omega$ is extremal if and only if its GNS representation $\pi_\omega$ is factorial, and
\item[(ii)] a full corner $pMp$ is a factor of type $t$ if and only if $M$ is a factor of type $t$. 
\end{enumerate}
Other conclusions are standard. \hfill $\Box$

A few remarks are in order.

  \begin{rmrk}
  \begin{enumerate}
  Let $\beta$ be a non-zero real number. 
  \item[(1)]
  Let $\omega$ be a $\beta$-KMS state. Let $m:=m_{\omega}$ be the $e^{-\beta c}$-conformal measure on $\overline{Y}_u$ associated to $\omega$ as in Prop. \ref{statetomeasure}. Then, $m$ is concentrated on $Y_u:=\overline{Y}_u \backslash \{G\}$. Note that, for every $s \in G$, 
 
 \[
 m(\{G\})=m(\{G\}+s)=e^{-\beta c(s)}m(\{G\}).\]
 Since, $\{G\}$ is of finite measure (as $\{G\} \subset \overline{X}_u$) and $\beta c$ is a non-zero homomorphism, we have $m(\{G\})=0$. 
  
  \item[(2)]  Suppose $\omega$ is  $\beta$-KMS state on $A:=C_{c}^{*}(P)$. Assume that $\omega$ is extremal and denote the GNS representation of $A$ associated to $\omega$ by $\pi_\omega$. Denote the associated cyclic vector by $\Omega_\omega$. Then, the factor $M:=\pi_\omega(A)^{''}$ is neither of type $I_n$ for $n$ finite nor of type $II_1$.
  
     To see this,  let $m$ be the $e^{-\beta c}$-conformal measure on $\overline{Y}_u$ associated to $\omega$. Pick $a \in P$ such that $c(a) \neq 0$. Note that $\pi_\omega(e_a)=\pi_\omega(1_{\overline{X}_u+a})$. Here $e_a=v_av_a^*$.  
  Then, \[\langle \pi_\omega(1-e_a)\Omega_\omega|\Omega_\omega \rangle=1-m(\overline{X}_u+a)=1-e^{-\beta c(a)} \neq 0.\] This implies that $\pi_\omega(1-e_a) \neq 0$.
     Thus, $\pi_\omega(e_a)$ is a proper subprojection of $1$ which is Murray von-Neumann equivalent to $1$ in $M$. Hence the conclusion.
          \end{enumerate}
      \end{rmrk}
  
  In the next proposition, we work out the extremal $\beta$-KMS states whose associated measure is supported on an orbit. 
  
  \begin{ppsn}
  \label{type I extremal states}

\begin{enumerate}
\item[(1)] Let $m$ be an $e^{-\beta c}$-conformal measure on $Y_u$ supported on an orbit $Orb(A)$ for some $A$. Denote the stabiliser of $A$ by $H$, i.e
 \[H:=\{s \in G: A+s=A\}.\] Let $\chi$ be a character of $H$. Define a state $\omega_{\chi,m}$ on $C_{c}^{*}(P)$ by
 \begin{equation*}
\omega_{\chi,m}(f w_s):=\begin{cases}
 0  & \mbox{ if
} s \notin H \cr
   &\cr
    \overline{\chi(s)}\displaystyle \int_{X_u} f(A)dm(A) &  \mbox{ if $s \in H$}.
         \end{cases}
\end{equation*}
Then, $\omega_{\chi,m}$ is an extremal $\beta$-KMS state on $C_{c}^{*}(P)$. Moreover, $\omega_{\chi,m}$ is of type I.
\item[(2)] Let $\omega$ be an extremal $\beta$-KMS state on $C_{c}^{*}(P)$ and let $m$ be the  $e^{-\beta c}$-conformal measure on $Y_u$ associated to $\omega$. Suppose that $m$ is atomic and concentrated on $Orb(A)$ for some $A \in Y_u$. 
Let $H$ be the stabiliser of $A$. Then, $\omega$ is of the form $\omega_{\chi,m}$ for some character $\chi$ of $H$. 

\end{enumerate}

\end{ppsn}
\textit{Proof.}
	 The fact that $\omega:=\omega_{\chi,m}$ is a $\beta$-KMS states follows directly from Theorem 1.3 of \cite{Neshveyev}. Let $\pi_\omega$ be  the GNS representation of
$\omega_{\chi,m}$. To prove that $\omega_{\chi,m}$ is extremal of type I, it suffices to show that $\pi_{\omega}(C^*(\overline{X}_u \rtimes P))''$ is a factor of type I. The proof  is similar to the proof of Thm. \ref{main_sec3}.

By the conformality of $m$, for a Borel set $E\subset Orb(A)$ with $m(E)>0$ and $t \in H$, we have $m(E)=m(E+t)=e^{-\beta c(t)}m(E)$ implying that $c(t)=0$. Thus, $H\subset c^{-1}(0)$.
 Let 
$\lambda=\{\lambda_t\}_{t \in G}$ be the Koopman representation of $G$ on $L^2(\overline{Y}_u)$. Note that since $Y_u=Orb(A)$ upto a null set, $L^2(\overline{Y}_u)=L^2(Orb(A))$ and  we can conclude that $\lambda_t=\lambda_s$, whenever $t-s \in H$.  Extend the character $\chi$ of $H$ to a character of $G$ which we denote again by $\chi$.
	
		Let $\clk$ be the Hilbert space defined by $\clk=l^2(G/H,L^2(\overline{Y}_u))$ and $\pi_0$ be the representation of $C_0(\overline{Y}_u)\rtimes G$ on $\clk$ given by
	\begin{eqnarray*}
		\pi_0(fu_t)(\xi)(\overline{s})=\overline{\chi(t)} f\lambda_{t}(\xi(\overline{s-t}))
	\end{eqnarray*}
Clearly, $\pi_0(C_0(\overline{Y}_u)\rtimes G)''=L^\infty(\overline{Y}_u)\rtimes G/H$.
 Let $\clh$ be the Hilbert subspace of $\clk$ given by \[\clh:=\{\xi \in l^2(G/H,L^2(\overline{Y}_u)): \xi(\overline{t})\in L^2(\overline{X}_u+t)\}.\] We have $\overline{s}=\overline{t}$ implies $\overline{X}_u+s=\overline{X}_u+t$ upto a null set, and hence the definition of $\clh$ makes sense. Since the crossed product $L^\infty(\overline{Y}_u)\rtimes G/H$ leaves $\clh$ invariant, 
we get a representation $\overline{\pi}: L^\infty(\overline{Y}_u)\rtimes G/H \to B(\clh)$. Define, $\pi=\overline{\pi}\circ\pi_0$. Consider the Hilbert space $P\clh$ where $P=\pi(1_{\overline{X}_u})$.  Let $\xi \in P\clh$ be given by,
\begin{equation*}
\xi(\overline{t})=\begin{cases}
1_{\overline{X}_u}  & \mbox{ if~~} \overline{t}=0,\\
0 &  \mbox{ otherwise }.
\end{cases}
\end{equation*} 
Let $\widetilde{\pi}: C^{*}(\overline{X}_u \rtimes P) \to B(P\clh)$ be defined by $\widetilde{\pi}(fw_t)=P\pi(\tilde{f}u_t)P$ where $\tilde{f}=f1_{\overline{X}_u}$. 

Arguing as before, we can conclude the following.
\begin{enumerate}
	\item The map $\overline{\pi}$ is injective and hence the von Neumann algebra  $\pi(C_0(\overline{Y}_u)\rtimes G)^{''}$ is isomorphic to $L^\infty(\overline{Y}_u)\rtimes G/H$.
	\item The triple $(\widetilde{\pi},P\clh,\xi)$ is the GNS representation of $\omega_{\chi,m}$.
	\item The von Neumann algebra generated by $\widetilde{\pi}(C^*(\overline{X}_u \rtimes P))$ is isomorphic to the full corner $1_{\overline{X}_u}
	(L^\infty(\overline{Y}_u)\rtimes G/H)1_{\overline{X}_u}$.
\end{enumerate}
Thus, $\pi_{\omega}(C^*(\overline{X}_u \rtimes P))''$ is isomorphic to the full corner $1_{\overline{X}_u}
(L^\infty(\overline{Y}_u)\rtimes G/H)1_{\overline{X}_u}$ which is clearly a factor of type I (as the measure $m$ is atomic). This completes the proof of $(1)$.

Conversely, suppose $\omega$ is an extremal $\beta$-KMS state and let $m$ be the $e^{-\beta c}$-conformal measure on $Y_u$ associated to $\omega$. Assume that $m$ is atomic and concentrated on $Orb(A)$ for some $A\in Y_u$. Let $H$ be the stabilizer of $A$. 
Appealing to Theorem 1.3 in \cite{Neshveyev} for the groupoid $\clg=\overline{X}_u \rtimes P$, we see that $\omega$ is determined by the measure $m$ and  a  measurable field of  states $\{\phi_B\}_{B \in \overline{X}_u}$ (each $\phi_B$ is a state on $C^{*}(\clg_B^B)$) satisfying certain compatibility conditions  such that 
\begin{equation}\label{neshveyev1}
\omega(f)=\int_{\clg^{(0)}}\sum_{\gamma \in \clg^B_B}f(\gamma)\phi_B(u_\gamma)~dm\end{equation}
for $f \in C_{c}(\clg)$.
Since $Y_u=Orb(A)$ upto a null set, the field of states $\{\phi_B\}_{B}$ reduces to a single state $\phi$ on $C^*(H)$. Now 
rewriting Eq. \ref{neshveyev1} for $fw_t$, we get
\begin{equation}\label{neshveyev2}
\omega(fw_t)=\begin{cases}
 \phi(u_{-t})\int_{X_u}f dm & \mbox{ if~~} t \in H,\\
0 &  \mbox{ otherwise}.
\end{cases}
\end{equation}
where $u_t, ~t \in H$ are the canonical unitaries of $C^*(H)$. 
Since $\omega$ is extremal, $\phi$ is extremal and consequently, $\phi$ is a character of $H$ which we denote by $\chi$. Then, $\omega=\omega_{\chi,m}$. Hence the proof. \hfill $\Box$

We end this section by explaining a recipe that allows us to construct conformal measures on the $(G,P)$-space $(Y_u,X_u)$, and consequently KMS states on $C_c^{*}(P)$ of the desired type.
 Let $(Y,X)$ be a pure $(G,P)$-space. For $y \in Y$, let \[Q_y:=\{s \in G: y-s \in X\}.\] First, let us check that for $y \in Y$, $Q_y \in Y_u$. Fix $y \in Y$. 
\begin{enumerate}
\item[(1)] Since, $\displaystyle Y=\bigcup_{a \in P}(X-a)$, it follows that there exists $a \in P$ such that $y \in X-a$. Then, $-a \in Q_y$. Thus, $Q_y$ is non-empty. 
\item[(2)]  As the intersection $\displaystyle \bigcap_{a \in P}(X+a)=\emptyset$, it follows that there exists $a \in P$ such that $y \notin X+a$. Then, $a \notin Q_y$ for such an element $a$. Hence, $Q_y$ is a proper subset of $G$. 
\item[(3)] The fact that $X+P \subset X$ implies that $-P+Q_y \subset Q_y$.
\end{enumerate}

It is clear that  the map 
$
Y \ni y \to Q_y \in Y_u$, denoted $T$, 
is  $G$-equivariant. 
  Note that for $s \in G$, 
\[
1_{Q_y}(s)=1_{X+s}(y).\]
This implies that the map $T$ is measurable.   Observe that  $T^{-1}(X_u)=X$. 

Let $(Y,X)$ be a pure $(G,P)$-space and let $m$ be an $e^{-\beta c}$-conformal measure on $Y$.  Clearly, the push-forward $T_*m:=m \circ T^{-1}$ is an $e^{-\beta c}$-conformal measure on the pure $(G,P)$-space $(Y_u,X_u)$. 
Suppose, in addition, that  $T$ is injective. Then, $T$ sets up a $G$-equivariant isomorphism between $Y$ and $T(Y)$. Consequently, up to null sets, we can identify, via the map $T$, the dynamical systems $(Y_u,T_*m,G)$ and $(Y,m,G)$. Moreover, $X_u$ gets identified with $X$. 

Thus, we have the following.  
\begin{rmrk}
\label{defining example}
Suppose $t  \in \{II,III\}$. The problem of constructing free and ergodic $e^{-\beta c}$-conformal measures on $Y_u$  of type $t$ is equivalent to the problem of exhibiting pure $(G,P)$-spaces $(Y,X)$, together with an   $e^{-\beta c}$-conformal measure $m$ on $(Y,X)$, such that 
\begin{enumerate}
\item[(A1)] the map
\[
Y \ni y \to Q_y:=\{s \in G: y-s \in X\} \in Y_u\]
is injective, and
\item[(A2)] the dynamical system $(Y,m,G)$ is free, ergodic and is of type $t$. 
\end{enumerate} 
It is also clear that metrically non-isomorphic dynamical systems give rise to distinct measures on $Y_u$. The notion of metric isomorphism for $(G,P)$-spaces are defined in the usual way as follows. 

For $i=1,2$, let $(Y_i,X_i)$ be a pure $(G,P)$-space and let $m_i$ be an $e^{-\beta c}$-conformal measure on $(Y_i,X_i)$. We say that $(Y_1,X_1,m_1)$ and $(Y_2,X_2,m_2)$ are metrically isomorphic, 
or simply isomorphic, if there exist $G$-invariant null sets $N_1 \subset Y_1$, $N_2 \subset Y_2$ and an invertible measurable map $S:Y_1\backslash N_1 \to Y_2 \backslash N_2$ such that 
\begin{enumerate}
\item[(i)] the map $S$ is $G$-equivariant, $S(X_1\backslash N_1)=X_2\backslash N_2$, and
\item[(ii)] for every Borel subset $E \subset Y_2 \backslash N_2$,  $m_2(E)=m_1(S^{-1}(E))$.
\end{enumerate}
\end{rmrk}

\section{The case $P=\mathbb{N}^{2}$}
In this section, we discuss the structure of KMS states on $C_{c}^{*}(P)$ for $\sigma^{c}$ when $P=\bbn^2$. Let us fix notation. Set $e_1:=(1,0)$ and $e_2:=(0,1)$. Define $v_1:=e_1$ and $v_2:=e_1+e_2$. 
Let $c:\bbz^2 \to \bbr$ be a non-zero homomorphism. We normalise and assume that $c(e_1)=1$ and $c(e_2)=\theta$. This does not lead to any loss of generality. In what follows, the homomorphism $c$ will 
be fixed.

First, we obtain a reasonable parametrisation of $Y_u$. Let $\Omega:=\{0,1\}^{\bbz}$ be the Cantor space.  Let $\tau:\Omega\to\Omega$ be the Bernoulli shift defined by \[\tau(x)_k:=x_{k-1}.\] 
Define a  $\bbz^2$-action on $\Omega\times \bbz$ by the following formulae.
\begin{align*}
 \big(x,t\big)+v_1&=\big(\tau(x),t+x_{-1}\big),\text{ and }\\
  \big(x,t\big)+v_2&=\big(x,t+1\big).
\end{align*}
For $(x,t)\in\Omega\times \bbz$, define $a(x,t):=a=(a_m)_{m\in\bbz}$ as follows.
\begin{equation}
 a_m:=\begin{cases}
  t-(x_0+x_1+\cdots +x_{m-1}) & \mbox{if~} m> 0, \cr
      t & \mbox{if $m=0$~}, \cr
    t+(x_{-1}+x_{-2}+\cdots +x_{m}) & \mbox{if~} m<0.
         \end{cases}
\end{equation}
Note that $a_m-a_{m+1}=x_m$. 
 Let $A(x,t)$ be defined by \[A(x,t):=\{mv_1+nv_2:  n \leq a_m\}.\] It is not difficult to verify that  $-e_1+A(x,t) \subset A(x,t)$ and $-e_2+A(x,t) \subset A(x,t)$. In other words, $A(x,t) \in Y_u$ for every $(x,t) \in \Omega \times \bbz$. 
 
 \begin{ppsn}
 \label{bijection}
 With the foregoing notation, the map \[\Omega\times \bbz\ni(x,t)\mapsto A(x,t)\in Y_u\] is a $\bbz^2$-equivariant homeomorphism.
  \end{ppsn}
\textit{Proof.} It is routine to check that the prescribed map is $\bbz^2$-equivariant. First, we show that it is onto.  
Let $A \in Y_u$ be given. By translating, if necessary, we can assume that $0 \in A$ which implies that $-\bbn^2 \subset A$. 

Let $m \in \bbz$ be given. We claim that  the set $\{k \in \bbz: mv_1+kv_2\in A\}$ is non-empty and bounded above. Choose $k_0<0$  such that $m+k_0<0$. Then, we have 
\[ mv_1+k_0v_2=(m+k_0)e_1+k_0e_2\in -\bbn^2\subseteq A.\] 
Hence, the set $\{k \in \bbz: mv_1+kv_2\in A\}$ is non-empty.
Suppose  $\{k \in \bbz: mv_1+kv_2\in A\}$ is not bounded above. Then, for every $k\in\bbn$, there exists a natural number $n_k\geq k$ such that $mv_1+n_kv_2\in A$.
Since $A-\bbn^2 \subset A$, we have $mv_1+(n_k-\ell)v_2\in A$ for every $\ell \in \bbn$ and for every $k \in \bbn$. This means that $mv_1+nv_2\in A$ for every $n\in\bbz$.
Again for any $k,\ell\in\bbn$ and $n\in \bbz$, we have
\[(m+n-k)e_1+(n-\ell)e_2=mv_1+nv_2-ke_1-\ell e_2\in A-\bbn^2 \subset A.\]
This implies that $A=\bbz^2$ which is contradiction to the fact that $A \in Y_u$. Therefore, the set $\{k \in \bbz: mv_1+kv_2\in A\}$ is bounded above. The proof of the claim is now complete. 

For $m\in \bbz$, define \[a_m:=\text{max}\{k \in \bbz: mv_1+kv_2\in A\}.\] Then, $a_m$ is an integer.  Let $A_a=\{mv_1+nv_2: n \leq a_m\}$. 
We claim that $A=A_a$ and for every $m \in \bbz$, $0 \leq a_{m}-a_{m+1} \leq 1$.

Clearly, $A \subset A_a$. Suppose  $mv_1+nv_2\in A_a$.  By the definition of $a_m$, $mv_1+a_mv_2\in A$. Note that
\[
 mv_1+nv_2 =mv_1+a_mv_2+(n-a_m)(e_1+e_2)  \in A-\bbn^2\subseteq A.\]
Hence, $A_a\subseteq A$. Therefore, $A=A_a$.

Let $mv_1+nv_2\in A_a$ be given. Since $A_a-\bbn^2 \subset A_{a}$, we have $mv_1+nv_2-ke_1-\ell e_2\in A_a$, for every $k,l\geq 0$.
Therefore, whenever $mv_1+nv_2\in A_a$, we have $a_{m-k+\ell}\geq n-\ell$.
In particular, when $n=a_m,k=1$ and $\ell=0$, we have $a_{m-1}\geq a_m$, and when $n=a_m,k=0$ and $\ell=1$, we have $a_{m+1}\geq a_m-1$.
Thus, $0\leq a_m-a_{m+1}\leq 1$ for every $m\in \bbz$.
The proof of the claim is now over. 

Define $x \in \Omega$ by setting $x_m:=a_{m}-a_{m+1}$ and let $t:=a_0$. Clearly, $A_a=A(x,t)$. Consequently, $A=A(x,t)$. This proves the surjectivity of the map
\[
\Omega \ni (x,t) \to A(x,t) \in Y_u.\]
The rest of the assertions require routine verifications which we leave to the reader. This completes the proof. \hfill $\Box$

\begin{rmrk}
\label{countably many non-trivial stabilisers}
For $A \in Y_u$, let $G_A$ be the stabiliser of $A$, i.e. \[G_A:=\{(m,n) \in \bbz^2: A+(m,n)=A\}.\] Let $(x,t) \in \Omega \times \bbz$ and let $A:=A(x,t) \in Y_u$ be the set defined as in 
Prop. \ref{bijection}. Then, $G_A \neq 0$ if and only if $x$ is a periodic point, i.e. there exists $p>0$ such that $x_{m+p}=x_m$ for all $m \in \bbz$. Moreover, there are only countably many periodic points in $\Omega$. This has the consequence that the groupoid 
$\mathcal{G}:=\overline{X}_u \rtimes \bbn^2$ (and also the transformation groupoid $\overline{Y}_u \rtimes \bbz^2$) has only countably many points in its unit space whose stabiliser is non-trivial. 
\end{rmrk}

We identify $Y_u$ with $\Omega \times \bbz$ via the map prescribed in Prop. \ref{bijection}. We abuse notation and write $Y_u=\Omega \times \bbz$. Then, $X_u=\Omega \times \bbn$. 
Let $m$ be a probability measure on $\Omega$.  Define a measure $\overline{m}$  on $\Omega \times \bbz$ by setting
\[\overline{m}(E\times \{n\})=(1-e^{-\beta(1+\theta)})e^{-\beta (1+\theta)n}m(E)\]
for a measurable subset $E \subset \Omega$. 

 Define $\chi:\Omega\to\bbr$ by \begin{equation*}
\chi(x):=\begin{cases}
 1  & \mbox{ if
} x_{-1} =0,\cr
   &\cr
   - \theta &  \mbox{ if } x_{-1}=1.
         \end{cases}
\end{equation*}
Let $\beta$ be a real number and let $m$ be a probability measure on $\Omega$. The measure $m$ is said to be $e^{-\beta \chi}$-conformal for the Bernoulli shift $\tau$ if for every Borel subset $E \subset \Omega$, 
\[
m(\tau(E))=\int_{E} e^{-\beta \chi}dm.\]

  \begin{ppsn}
  \label{bijection 2}
Suppose $\beta (\theta+1)>0$. Then, the map $m \to \overline{m}$ defines a bijection between the set of $e^{-\beta \chi}$-conformal measures on $\Omega$ and
the set of $e^{-\beta c}$-conformal measures on the $(\bbz^2,\bbn^2)$-space $(Y_u,X_u)=(\Omega \times \bbz,\Omega \times \bbn)$. 
\end{ppsn}
\textit{Proof.} 
The proof is not difficult. We have included some details for completeness. Suppose that $m$ is an $e^{-\beta \chi}$-conformal measure on $\Omega$. Clearly, $\overline{m}(X_u)=1$. 
Let $F \subset \Omega$ be  a Borel subset and let $n \in \bbz$ be given. It suffices  to show that
\begin{align*}
 \overline{m}((F\times \{n\})+v_1)&=e^{-\beta }\overline{m}(F\times \{n\}),~\textrm{ and }\\
  \overline{m}((F\times \{n\})+v_2)&=e^{-\beta (1+\theta)}\overline{m}(F\times \{n\}).
\end{align*}
Define $F_0:=\{x \in F: x_{-1}=0\}$ and $F_1:=\{x \in F: x_{-1}=1\}$. 
Calculate as follows to observe that 
\begin{align*}
 \frac{1}{1-e^{-\beta(1+\theta)}}\overline{m}((F\times \{n\})+v_1)
 &=\frac{1}{1-e^{-\beta(1+\theta)}}\Big(\overline{m}(\tau(F_0)\times \{n\})+\overline{m}(\tau(F_1)\times \{n+1\})\Big)\\
  &=e^{-\beta n(1+\theta)}  m(\tau(F_0))+ e^{-\beta (n+1)(1+\theta)}m(\tau(F_1))\\
  &=e^{-\beta n(1+\theta)} e^{-\beta} m(F_0)+ e^{-\beta (n+1)(1+\theta)} e^{\beta \theta}m(F_1)\\
  &=e^{-\beta n(1+\theta)} e^{-\beta}\big(m(F_0)+ m(F_1)\big)\\
  &=e^{-\beta}\frac{1}{1-e^{-\beta(1+\theta)}}\overline{m}(F\times \{n\}).
\end{align*}
Similarly, we can prove that $ \overline{m}((F\times \{n\})+v_2)=e^{-\beta (1+\theta)}\overline{m}(F\times \{n\})$.
Hence,  $\overline{m}$ is an $e^{-\beta c}$-conformal measurs on $(\Omega\times \bbz,\Omega \times \bbn)$.

Conversely, suppose $\mu$ is an $e^{-\beta c}$-conformal measure on $(Y_u,X_u)$. Define a measure $m$ on $\Omega$ by 
\[
m(E)=\frac{1}{1-e^{-\beta(1+\theta)}}\mu(E \times \{0\})\]
for a Borel subset $E \subset \Omega$. Using the conformality condition on $\mu$ and the fact that $E \times \{n\}=(E \times \{0\})+nv_2$, it is routine to see that $\overline{m}=\mu$. 

Let $F \subset \Omega$ be measurable.  
Set $F_0:=\{x\in F: x_{-1}=0\}$ and $F_1:=\{x\in F: x_{-1}=1\}$. Thanks, to the conformality condition on $\mu$ and the calculation done earlier, it follows that the equality
$
\mu((F \times \{0\})+v_1)=e^{-\beta}\mu(F \times \{0\})$ is equivalent to the equality 
\[  m(\tau(F_0))+ e^{-\beta(1+\theta)}m(\tau(F_1))=   e^{-\beta}\big(m(F_0)+ m(F_1)\big)=e^{-\beta}m(F).\]
In other words, for every Borel subset $F \subset \Omega$, 
\[
m(F)=\int_{\tau(F)}e^{\beta (\chi \circ \tau^{-1})}dm.\]
Hence, $\frac{d(m\circ \tau)}{dm}=e^{-\beta \chi}$, i.e. $m$ is $e^{-\beta \chi}$-conformal. The fact that the map $m \to \overline{m}$ is a bijection is clear.  \hfill $\Box$

\begin{rmrk}
\label{necessity condition}
 For the existence of an $e^{-\beta c}$-conformal measure on $(Y_u,X_u)$, it is necessary that $\beta(1+\theta)>0$. This is because, if $\overline{m}$ is an $e^{-\beta c}$-conformal measure on $Y_u$, then
\[
1=\overline{m}(X_u)=\sum_{n=0}^{\infty}\overline{m}(\Omega \times \{n\})=\sum_{n=0}^{\infty}\overline{m}((\Omega \times \{0\})+nv_2)=\sum_{n=0}^{\infty}e^{-\beta(1+\theta)n}\overline{m}(\Omega \times \{0\}).\]
 Hence, the condition $\beta(1+\theta)>0$ is necessary. 

 The bijection $m \to \overline{m}$ of Prop. \ref{bijection 2} preserves freeness, ergodicity and type. Thus,  $(\Omega,\tau,m)$ is free
if and only if $(Y_u,\bbz^2,\overline{m})$ is free. Similarly, $(\Omega,\tau,m)$ is ergodic and is of type $t$ if and only if $(Y_u,\bbz^2,\overline{m})$ is ergodic and is of type $t$. 
Here, $t \in \{I,II,III\}$. 

\end{rmrk}

We are now in a position to determine the values of $\beta$ and $\theta$ for which there is a $\beta$-KMS state on $C_{c}^{*}(\bbn^2)$ for $\sigma^{c}:=\{\sigma_t\}_{t \in \bbr}$. Recall that the flow
$\sigma^{c}:=\{\sigma_t\}_{t \in \bbr}$ on $C_{c}^{*}(\bbn^2)$ is defined by 
\[
\sigma_t(v_{(m,n)})=e^{i(m+n\theta)t}v_{(m,n)}.\]

\begin{ppsn}
\label{values}
Suppose $\beta$ is a non-zero real number. The following are equivalent. 
\begin{enumerate}
\item[(1)] There is a $\beta$-KMS state on $C_{c}^{*}(\bbn^2)$ for $\sigma^c$. 
\item[(2)] There is an $e^{-\beta c}$-conformal measure on $(Y_u,X_u)$. 
\item[(3)] There is an $e^{-\beta \chi}$-conformal measure on $\Omega$ and $\beta(1+\theta)>0$.
\item[(4)] $\beta>0$ and $\theta\geq 0$. 
\end{enumerate}
\end{ppsn}
\textit{Proof.} The equivalence between $(1)$, $(2)$ and $(3)$ follow from  Prop. \ref{measuretostate}, Prop. \ref{bijection 2} and Remark \ref{necessity condition}. We prove the equivalence between $(3)$ and $(4)$. 
Assume that $(3)$ holds. Let $m$ be an $e^{-\beta \chi}$-conformal measure. 
Suppose $\theta<0$. Then, the potential $\beta \chi$ is strictly negative or strictly positive. Suppose $\beta \chi$ is strictly positive. Then, there exists a real number $a>0$ such that $\beta \chi>a$. Note that
\[
1=m(\tau(\Omega))=\int_{\Omega}e^{-\beta \chi}dm \leq e^{- a}<1.\]
This is a contradiction. A similar contradiction will be met if $\beta \chi$ is strictly negative. Therefore, $\theta \geq 0$. (A direct application of Thm. 6.2 of \cite{Thomsen} could also have been made instead). The condition $\beta(1+\theta)>0$ implies that $\beta>0$. 
This completes the proof of $(3) \implies (4)$. 
  
  Assume now that $\beta>0$ and $\theta \geq 0$. Let us fix some notation. For $n \in \bbz$, define $S_{n}(\chi)$ by setting
\begin{equation*}
S_n(\chi)(x):=\begin{cases}
 \sum_{j=0}^{n-1} \chi \circ \tau^{j}(x)  & \mbox{ if
} n \geq 1,\cr
0 & \mbox{ if } n=0 \cr
   - \sum_{j=1}^{|n|}\chi \circ \tau^{-j}(x) &  \mbox{ if } n \leq -1.
         \end{cases}
\end{equation*}
Note that if $n \leq -1$, then $S_{n}(\chi)(x)=-|n|+(1+\theta)(x_{0}+x_{1}+\cdots+x_{|n|-1})$, and if $n \geq 1$, we have $S_n(\chi)(x)=n-(1+\theta)(x_{-1}+x_{-2}+\cdots+x_{-n})$. 

We apply Lemma 4.9 of \cite{Thomsen} which states the following. Let $x\in\Omega$ be given. If $x$ is periodic of period $p>0$, then there is an $e^{-\beta \chi}$-conformal measure on $\Omega$ concentrated on the orbit of $x$ iff $S_p(\chi)(x)=0$. If $x$ is not periodic, then  there is an $e^{-\beta  \chi}$-conformal measure concentrated on the orbit of $x$ if and only if 
 $\displaystyle \sum_{n\in\bbz} e^{-\beta S_n(\chi)(x)}<\infty$.
 
 \textbf{Case 1:} Suppose $\theta=0$. Let $x\in \Omega$ be such that $x_m=1$ for every $m \in \bbz$. Clearly, $x$ is periodic of period $1$ and $S_{1}(\chi)(x)=0$. Therefore, there is an $e^{-\beta \chi}$-conformal measure supported on the orbit of $x$. 
 
 \textbf{Case 2:} Suppose $\theta>0$. We exhibit uncountably many non-periodic points $x$ for which the series $\displaystyle \sum_{n \in \bbz}e^{-\beta S_n(\chi)(x)}$ converges. Let $\ell$ be a positive integer such that $\ell>1+\theta$.  Suppose   $S \subset \ell \bbn$. Denote the indicator function of $S$ by $1_{S}$.  Define $x \in \Omega$ by setting 
 \begin{equation*}
 x_m:=\begin{cases}
  1 & \mbox{ if~} m \geq 0 \cr
   &\cr
    1_{S}(-m) & \mbox{if~} m<0
         \end{cases}
\end{equation*}
Note that $S_{n}(\chi)(x)=|n|\theta$ if $n \leq -1$. Therefore, the series $\sum_{n\leq -1}^{\infty}e^{-\beta S_n(\chi)(x)}$ is convergent. 

For $n \geq 1$, observe that since $x_{-1}+x_{-2}+\cdots+x_{-n} =|S \cap \{1,2,\cdots,n\}| \leq \frac{n}{\ell}$,
\[
S_{n}(\chi)(x)=n-(1+\theta)(x_{-1}+x_{-2}+\cdots+x_{-n}) \geq n(1-\frac{1+\theta}{\ell}).\] 
Thus, $
e^{-\beta S_n(\chi)(x)} \leq e^{-\beta n (1-\frac{1+\theta}{\ell})}$. Consequently, the series $\sum_{n \geq 1}^{\infty}e^{-\beta S_n(\chi)(x)}$ converges. Thanks to Lemma 4.9 of \cite{Thomsen}, there is an $e^{-\beta \chi}$-conformal measure supported on the orbit of $x$. The proof is complete. \hfill $\Box$

We discuss the structure of KMS states in more detail now. Assume hereafter that $\beta$ is an arbitrary positive real number. 

\textbf{The case $\theta=0$:} Assume that $\theta=0$. Denote the element in $\{0,1\}^{\bbz}$ whose every entry is 1 by $\underline{1}$. Suppose that $m$ is an $e^{-\beta \chi}$-conformal measure on $\Omega$.  We claim that $m$ is concentrated on $\underline{1}$. 
 We claim that  $m(\{x:x_{-1}=0\})=0$. Suppose  $m(\{x:x_{-1}=0\})>0$. Then, by conformality, we have
\begin{align*}
1&=m(\tau(\Omega))\\
&=\int e^{-\beta \chi}dm\\
&=e^{-\beta}m(\{x:x_{-1}=0\})+m(\{x:x_{-1}=1\})\\ 
&<m(\{x:x_{-1}=0\})+m(\{x:x_{-1}=1\})=1\end{align*}
which is a contradiction. 

Therefore, for almost all $x$, $x_{-1}=1$. Since $m$ is quasi-invariant for the Bernoulli shift $\tau$, we have, for every $k$, $x_k=1$ for almost all $x$. Hence, $x=\underline{1}$ for almost all $x$. This proves the claim. 

Let $\omega$ be a $\beta$-KMS state on $C_{c}^{*}(\bbn^2)$. Thanks to Prop. \ref{bijection} and Prop. \ref{bijection 2} and by what we have proved now, its associated $e^{ - \beta c}$-conformal measure $m_{\omega}$
is supported on $Orb(A(\underline{1},0))$. Note that \[A:=A(\underline{1},0)=\{(m,n): m \geq 0\}\] whose stabiliser $G_{A}=\{0\}\times \bbz=c^{-1}(\{0\})$. Appealing to Theorem 1.3 of \cite{Neshveyev}, we see that there
exists a state $\phi$ on $C^{*}(G_A)$, or equivalently a probability measure $\mu$ on $\bbt$ such that 
\begin{equation}
\label{KMS theta zero}
\omega(fw_{(m,n)})=\delta_{m,0} \Big(\int z^{n}d\mu(z)\Big) (1-e^{-\beta})\sum_{k=0}^{\infty}e^{-\beta k}f(A(\underline{1},k)).
\end{equation}
Conversely, given a probability measure $\mu$ on $\bbt$, if we define $\omega$ as in Eq. \ref{KMS theta zero}, then $\omega$ will be a $\beta$-KMS state. Thus, when $\theta=0$, the simplex of $\beta$-KMS states for $\sigma^{c}$ is homeomorphic 
to the simplex of probability measures on the circle. 

Hereafter, we assume that $\beta>0$ and $\theta>0$. 
We    summarise the   structure of extremal $\beta$-KMS states on $C_{c}^{*}(\bbn^2)$, that result from our discussions so far, as follows. 
\begin{enumerate}
\item[(1)] Suppose $\omega$ is an extremal $\beta$-KMS state on $C_{c}^{*}(\bbn^2)$. Denote the associated conformal measure on $Y_u$ by $\overline{m}$ and let $m$ be the corresponding $e^{-\beta \chi}$-conformal measure on $\Omega$ obtained via the map in Prop. \ref{bijection 2}. 
Then, thanks to Prop. \ref{ergodicity}, $m$ is ergodic. Suppose that $m$ is atomic and $m$ is concentrated on $Orbit(x)$ for some $x \in \{0,1\}^{\bbz}$. Set $A:=A(x,0)$.  Thanks to Prop. \ref{type I extremal states}, if the stabiliser $G_A=\{0\}$, then $\omega=\omega_{\overline{m}}$. 
              If $G_A \neq 0$, then $\omega$ is as in Prop. \ref{type I extremal states}. In both cases, $\omega$ is of type I.

 If $m$ is non-atomic, then $\omega=\omega_{\overline{m}}$ where $\omega_{\overline{m}}$ is the state corresponding to $\overline{m}$ obtained via the conditional expectation.
 This follows by applying Corollary 1.2 of \cite{Neshveyev} and by Remark \ref{countably many non-trivial stabilisers}.
               In this case, $\omega$ is of type II or III depending upon whether $m$ (or equivalently $\overline{m}$) is of type II or of type III.

 \item[(2)] Conversely, suppose  $m$ is an ergodic $e^{-\beta \chi}$-conformal measure on $\Omega$. Let $\overline{m}$ be the corresponding $e^{-\beta c}$-conformal measure on $Y_u$ given by Prop. \ref{bijection 2}. If $m$ is non-atomic, then $\omega_{\overline{m}}$ is an extremal $\beta$-KMS state and 
  its type is the same as that of $m$.    
  To see that $\omega_{\overline{m}}$ is extremal, it suffices to show that the GNS representation is factorial.  Thanks to Thm. \ref{main_sec3}, it suffices to prove that the $\bbz^2$-action on $Y_u$ is free. 
  
  Let $(m,n) \neq (0,0)$ be given. 
  Then, by Remark \ref{countably many non-trivial stabilisers}, the  set \[\{A \in Y_u :  A+(m,n) =A\}\] is countable. Since $\overline{m}$ is non-atomic, it follows that $\{A \in Y_u: A+(m,n)=A\}$, which is countable, is 
  a set of measure zero. Thus, the $\bbz^2$-action on $Y_u$ is free.

  Suppose $\overline{m}$ is atomic and $\overline{m}$ is concentrated on $Orbit(A)$ for some $A \in Y_u$. From Prop. \ref{type I extremal states}, if the stabiliser $G_A=\{0\}$, then $\omega_{\overline{m}}$ is an extremal $\beta$-KMS state of type I.
   If $G_A\neq \{0\}$, then $\omega_{\chi,m}$, as defined in Prop. \ref{type I extremal states}, is an extremal $\beta$-KMS state of type I.     
   
   \item[(3)] We have exhibited, in the proof of Prop. \ref{values}, uncountably many atomic probability measures on $\Omega$ that are $e^{-\beta \chi}$-conformal. Thus, there are uncountably many type I $\beta$-KMS states on $C_{c}^{*}(\bbn^2)$ for $\sigma^c$. 
\end{enumerate}

 We end this section by exhibiting a type II KMS state  when $\theta$ is irrational and a type III KMS state when $\theta=1$. We construct $(\bbz^2,\bbn^2)$-spaces with conformal measures satisfying $(A1)$ and $(A2)$ of Remark \ref{defining example}. Uncountably many such examples for every $\theta>0$ and for every $\beta>0$ will be constructed in Section 5. 

\textbf{A type II example:} Assume that $\theta \in (0,\infty)$ is irrational. Suppose $\beta>0$.  Let $Y:=\bbr$ and $X:=[0,\infty)$. Define a $\bbz^2$-action on $Y$ by 
\[
t+e_1:=t+1;~\textrm{and~} t+e_2:=t+\theta.\]
Clearly, $X+\bbn^2 \subset X$. Also, $(Y,X)$ is a pure $(\bbz^2,\bbn^2)$-space. Since $\theta$ is irrational, the $\bbz^2$-action on $Y$ is free. 
Let $m$ be the measure on $\bbr$ such that $dm:=\beta e^{-\beta t}dt$. Then, $m(X)=1$ and $m$ is $e^{-\beta c}$-conformal. As $\bbz+\bbz \theta$ is a dense subgroup of $\bbr$, 
the $\bbz^2$-action on $Y$ is ergodic. Clearly, it is of type II as $m$ is absolutely continuous w.r.t. the Lebesgue measure.

 Observe that, for $t \in \bbr$, 
\[
Q_t:=\{(m,n) \in \bbz^2: t-(me_1+ne_2) \in X\}=\{(m,n) \in \bbz^2: t \geq m+n\theta\}.\]
If $t_1<t_2$, the density of $\bbz+\bbz\theta$ in $ \bbr$ implies that there exists $(m,n) \in \bbz^2$ such that $t_1<m+n\theta<t_2$. Then, $(m,n) \in Q_{t_2}$ but $(m,n) \notin Q_{t_1}$. 
Therefore, the map \[Y \ni t \to Q_t \in Y_u\] is injective. 

Applying Remark \ref{defining example}, we get an extremal $\beta$-KMS state of type II when $\theta$ is irrational.

\textbf{A type III example:} Assume that $\theta=1$.  We make use of Arnold's dyadic adding machine to produce a type III example in this situation. Let us recall the basics on  adding machine from \cite{Aaronson}. Let $\beta>0$ be fixed. 
Let $p \in (0,\frac{1}{2})$ be such that $\frac{1-p}{p}=e^{\beta}$. 

Let \[
\Omega:=\prod_{n=1}^{\infty}\{0,1\}:=\{(x_1,x_2,\cdots,):x_n \in \{0,1\}\}\] be the group of dyadic integers. Let $\displaystyle \mu:=\otimes_{k=1}^{\infty}\mu_k$ be the product measure where the measure $\mu_k$ on $\{0,1\}$ is 
given by 
\[
\mu_k(\{0\})=1-p;~\mu_k(\{1\})=p.\]
We remove from $\Omega$ the null set of eventually constant sequences and we denote the resulting set again by $\Omega$. 

 Denote the map on $\Omega$ that corresponds to  addition by $\underline{1}$ by $\tau$. Recall that 
\[
\tau(1,1,1,\cdots,1,0,*,*,\cdots)=(0,0,0,\cdots,0,1,*,*,\cdots).\]

Define $\phi:\Omega \to \bbz$ by 
$
\phi(x):=\min\{n \geq 1: x_n=0\}-2$. Then, 
\[
\frac{d(\mu \circ \tau)}{d\mu}=e^{\beta \phi}.\]

Let $Y:=\Omega \times \bbz$ and $X:=\Omega \times \{0,1,2,\cdots\}$. Let $\overline{\mu}$ be the measure on $Y$ given by 
\[
d\overline{\mu}:=(1-e^{-\beta})e^{-\beta n}d\mu dn\]
where $dn$ is the counting measure on $\bbz$. 

Define a $\bbz^2$-action on $Y$ by 
\[
(x,t)+e_1:=(\tau x, \phi(x)+t+1);~\textrm{and~~} (x,t)+e_2:=(x,t+1).\]

Since $\phi \geq -1$, it follows that $X+\bbn^2 \subset X$. Also, it is routine to verify that $(Y,X)$ is a pure $(\bbz^2,\bbn^2)$-space and $\overline{\mu}$ is an $e^{-\beta c}$-conformal measure. 

\begin{lmma}
\label{separation}
Let $x,y \in \Omega$ be such that $x \neq y$. Then, there exists an  integer $m$ such that $\phi(\tau^{m}(x))\neq \phi(\tau^{m}(y))$. 
\end{lmma}
\textit{Proof.} Let $k$ be the least integer for which $x_k \neq y_k$. We can, without loss of generality, assume that $x_k=0$ and $y_k=1$. 
Let $m:=2^{k-1}-1-(x_1+2x_2+\cdots+2^{k-2}x_{k-1})$. A  calculation with `binary arithmetic' shows that 
\begin{align*}
\tau^{m}(x)&=(\underbrace{1,1,\cdots,1}_{k-1},0,*,*,*,\cdots) \\
\tau^{m}(y)&=(\underbrace{1,1,1,\cdots 1}_{k-1},1,*,*,*,\cdots).
\end{align*}
Thus, $\phi(\tau^m(x))=k-2 < \phi(\tau^m(y))$. The proof is complete. \hfill $\Box$

\begin{ppsn}
\label{adding machine}
With the foregoing notation, we have the following. 
\begin{enumerate}
\item[(1)] The $\bbz^2$-action on $Y$ is free and ergodic.
\item[(2)] The measure $m$ is of type III. 
\item[(3)] The map 
\[
Y \ni (x,t) \to Q_{(x,t)} \in Y_u\]
is injective. 
\end{enumerate}
\end{ppsn}
\textit{Proof.} The first two statements are straightforward consequences of Prop. 1.2.8 and Thm. 1.2.9 of \cite{Aaronson}. We include a brief explanation. 

Note that the $e_2$-action on $Y$ is by translation
by $1$ on the second coordinate. This forces that any $\bbz^2$-invariant subset of $Y$ must be of the form $E \times \bbz$ for some subset $E$ of $\Omega$.  The fact that $E \times \bbz$ is invariant
under the action of $e_1$ implies that $\tau(E)=E$. By Prop. 1.2.8 of \cite{Aaronson}, it follows that $E$ is either null or co-null. Hence, the $\bbz^2$-action is ergodic. Using the fact that the action of $\tau$ on $\Omega$ is free,
it is routine to check that the $\bbz^2$-action is free. 

Suppose $\overline{\mu}$ is not of type III. Let $\overline{\nu}$ be a $\sigma$-finite measure that is absolutely continuous w.r.t. $m$ such that $\overline{\nu}$ is $\bbz^2$-invariant. The invariance of $\overline{\nu}$ under $e_2$ implies that
$d\overline{\nu}:=d\nu dn$ for some $\sigma$-finite measure $\nu$ on $\Omega$ which is absolutely w.r.t. to $\mu$. The invariance of $\overline{\nu}$ under $e_1$ implies that $\nu$ is $\tau$-invariant which is a contradiction to
Theorem 1.2.9 of \cite{Aaronson}. The proof of $(1)$ and $(2)$ is complete. 

For $m \in \bbz$ and $x \in \Omega$, let $c(m,x) \in \bbz$ be defined by 
\begin{equation}
c(m,x):=\begin{cases}
 \sum_{k=0}^{m-1}\phi(\tau^k(x))  & \mbox{if
} m \geq 1,\cr
    0  & \mbox{if } m=0, \cr
    -\sum_{k=1}^{|m|}\phi(\tau^{-k}x) &  \mbox{if } m \leq -1.
         \end{cases}
\end{equation}

Let $(x,s), (y,t) \in Y$ be such that $Q_{(x,s)}=Q_{(y,t)}$. Notice that 
\[
Q_{(x,s)}=\{(m,n) \in \bbz^2:(x,s)-(me_1+ne_2) \in X\}=\{(m,n) \in \bbz^2:  n \leq -m+c(-m,x)+s\}.\]
The equality $Q_{(x,s)}=Q_{(y,t)}$ implies  in particular that for every $m \in \bbz$, \[-m+c(-m,x)+s=-m+c(-m,y)+t.\]
Substituting, $m=0$ in the above equality, we deduce that $s=t$. The fact that $c(-m,x)=c(-m,y)$ for every $m \in \bbz$ implies that $\phi(\tau^mx)=\phi(\tau^my)$ for each $m \in \bbz$. By Lemma \ref{separation}, we have
$x=y$. This proves $(3)$ and the proof is complete. \hfill $\Box$

Applying Remark \ref{defining example}, we see that when $\theta=1$, there is an extremal $\beta$-KMS state of type III.

\section{Examples of type II and type III}

In this section, we construct uncountably many  extremal $\beta$-KMS states of type II and of type III for every $\beta>0$ and for every $\theta>0$.  It suffices to exhibit uncountably many non-atomic ergodic $e^{-\beta c}$-conformal measures on the $(\bbz^2,\bbn^2)$-space $(Y_u,X_u)$ of the desired type. We do this by appealing to Remark \ref{defining example}. For the rest of this section, we assume 
that $\beta$ is an arbitrary positive real number.  Recall that $e_1=(1,0)$, $e_2=(0,1)$, $v_1=e_1$ and $v_2=e_1+e_2$. 

\textbf{Type II examples for an irrational $\theta$:} Assume that $\theta>0$ is irrational. Choose $\alpha \in (0,1)$ such that $1,\alpha$ and $\theta$ are rationally independent. For $\delta \in (0,\theta]$, let $C^{\delta}$ be the cone in $\bbr^2$ generated by $(0,1)$ and $(\alpha,\delta)$, i.e. 
\[
C^{\delta}:=\{x(0,1)+y(\alpha,\delta): x,y \geq 0\}.\]

Let $\pi:\bbr^2 \to \frac{ \bbr^2}{\bbz \times \{0\}}=\bbt \times \bbz$ be the quotient map. We denote $\pi(x,y)$ by $\overline{(x,y)}$. Denote the subgroup generated by $\overline{(\alpha,\theta)}$ and $\overline{(0,1)}$ by $\Gamma$. Let $\phi:\bbz^2 \to \Gamma$ be the homomorphism defined by \[\phi(e_1)=\overline{(0,1)}; ~~\phi(e_2)=\overline{(\alpha,\theta)}.\] Note that $\phi$ is an isomorphism as $\alpha$ and $\theta$ are irrational. 

Set $Y^{\delta}:=\frac{\bbr^2}{\bbz \times \{0\}}=\bbt \times \bbr$ and $X^{\delta}:=\pi(C^{\delta})$. It is not difficult to see that $X^{\delta}$ is closed. Define a $\bbz^2$-action on $Y^{\delta}$ by 
\begin{align*}
\overline{(x,y)}+e_1:&=\overline{(x,y)}+\phi(e_1)=\overline{(x,y+1)}\\
\overline{(x,y)}+e_2:&=\overline{(x,y)}+\phi(e_2)=\overline{(x+\alpha,y+\theta)}.
\end{align*}

Since, $\delta \in (0,\theta]$, it follows that $(\alpha,\theta) \in C^{\delta}$. The fact that $C^{\delta}$ is a cone implies that $X^{\delta}$ is invariant under $\bbn^2$, i.e. $X^{\delta}+\bbn^2 \subset X^{\delta}$. Thus, $(Y^{\delta},X^{\delta})$ is a $(\bbz^2,\bbn^2)$-space. We leave it to the reader to verify that it is pure. 

Define a measure $\mu$ on $Y^{\delta}$ by $d\mu:=e^{-\beta y}dx dy$.  Here, $dx$ is the Haar measure on $\bbt$ and $dy$ is the Lebesgue measure on $\bbr$. Note that $X^{\delta} \subset \bbt \times [0,\infty)$ and the latter set has finite $\mu$-measure. Thus, $\mu(X^{\delta})<\infty$. Define a measure $m$ on $Y^{\delta}$ by $dm:=\frac{1}{\mu(X^{\delta})}d\mu$. It is clear that $m$ is an $e^{-\beta c}$-conformal measure. 

\begin{rmrk}
\label{geometric}
We need the following `pictorially' obvious facts. We omit the rigorous proofs as they are elementary. 
\begin{enumerate}
\item[(1)] The interior $Int(X^\delta)$ is dense in $X^{\delta}$. 
\item[(2)] If $X^{\delta}+\overline{(x,y)}=X^{\delta}$ for some $\overline{(x,y)} \in \bbt \times \bbz$, then $\overline{(x,y)}=\overline{(0,0)}$. 
\item[(3)] Let $\delta_1,\delta_2 \in (0,\theta]$. Suppose there exists $\overline{(x,y)} \in \bbt \times \bbz$ for which 
$X^{\delta_1}+\overline{(x,y)}=X^{\delta_2}$. Then, $\delta_1=\delta_2$. 
\end{enumerate}
\end{rmrk}

\begin{ppsn}
Keep the foregoing notation. 
\begin{enumerate}
\item[(1)] The $\bbz^{2}$-action on $Y^{\delta}$ is free and ergodic. 
\item[(2)] The measure $m$ is of type II. 
\item[(3)] The map, denoted $T^{\delta}$,
\[
Y^{\delta} \ni \overline{(x,y)} \to Q_{\overline{(x,y)}} \in Y_u\]
is injective. 
\end{enumerate}
For $\delta \in (0,\theta]$, let $\mu^{\delta}:=T^{\delta}_{*}m$ be the push-forward measure on $Y_u$. If $\delta_1 \neq \delta_2$, then  $\mu^{\delta_1} \neq \mu^{\delta_2}$.
\end{ppsn}
\textit{Proof.} Note that the measure $m$ is absolutely continuous with respect to $dxdy$, which is the Haar measure on the cylinder $\bbt \times \bbz$, and $dxdy$ is clearly $\bbz^2$-invariant. As $1,\alpha$ and $\theta$ are rationally independent, the subgroup $\Gamma$ is dense in $\bbt \times \bbz$. Hence, the action of $\bbz^2$ on $Y^{\delta}$, which is by translations by $\Gamma$ via the isomorphism $\phi:\bbz^2 \to \Gamma$, is ergodic. As the map $\phi$ is an isomorphism, it is also clear that the $\bbz^2$-action on $Y^{\delta}$ is free. The proof of $(1)$ and $(2)$ are now complete. 

Let $(x_1,y_1), (x_2,y_2) \in \bbr^2$ be such that $Q_{\overline{(x_1,y_1)}}=Q_{\overline{(x_2,y_2)}}$. Note that for $i=1,2$, 
\begin{align*}
Q_{\overline{(x_i,y_i)}}:&=\{(m,n): \overline{(x_i,y_i)}-m\overline{(\alpha,\theta)}-n\overline{(0,1)} \in X^{\delta}\}\\
&=\{(m,n):m\overline{(\alpha,\theta)}+n\overline{(0,1)} \in -X^{\delta}+\overline{(x_i,y_i)}\}.
\end{align*}
The above equation, together with the equality $Q_{(\overline{x_1,y_1})}=Q_{\overline{(x_2,y_2)}}$ implies that 
\[
(-X^{\delta}+\overline{(x_1,y_1)}) \cap \Gamma= (-X^{\delta}+\overline{(x_2,y_2)}) \cap \Gamma.\]
Taking closure in the above equality, and using the fact that $\overline{Int(X^{\delta})}=X^{\delta}$ and $\Gamma$ is dense in $Y^{\delta}$, we deduce that 
$-X^{\delta}+\overline{(x_1,y_1)}=-X^{\delta}+\overline{(x_2,y_2)}$. By $(2)$ of Remark \ref{geometric}, we have $\overline{(x_1,y_1)}=\overline{(x_2,y_2)}$. This proves $(3)$. 

Let $\delta_1,\delta_2 \in (0,\theta]$. Suppose $\mu^{\delta_1}=\mu^{\delta_2}=\mu$.  Note that, for $i=1,2$, $\mu^{\delta_i}$ is concentrated on the image of $T^{\delta_i}$. 
Thus, up to a set of $\mu$-measure zero, $T^{\delta_1}(Y^{\delta_1})=T^{\delta_2}(Y^{\delta_1})$. Choose $(x_1,y_1), (x_2,y_2) \in \bbr^{2}$ such that 
$T^{\delta_1}(\overline{(x_1,y_1)})=T^{\delta_2}(\overline{(x_2,y_2)})$. This implies that 
\[(-X^{\delta_1}+\overline{(x_1,y_1)}) \cap \Gamma=(-X^{\delta_2}+\overline{(x_2,y_2)}) \cap \Gamma.\] Taking closure, we deduce 
\[
-X^{\delta_1}+\overline{(x_1,y_1)}=-X^{\delta_2}+\overline{(x_2,y_2)}.\] By $(3)$ of Remark \ref{geometric}, we have $\delta_1=\delta_2$. This completes the proof. \hfill $\Box$

\textbf{Type II examples for $\theta=1$:} Assume that $\theta=1$. Recall that the homomorphism $c:\bbz^2 \to \bbr$ is defined by $c(e_1)=1$ and $c(e_2)=1$. For an irrational $\alpha \in (0,1)$, let $R_\alpha$ be the rotation by angle $\alpha$ on $\bbt:=[0,1)$. 
In additive notation, $R_\alpha(x)=x+\alpha \mod 1$. For $\eta \in (0,1)$, let $h_\eta:=1_{[0,\eta)}$ and let $\phi_{\eta,\alpha}:=h_\eta \circ R_\alpha- h_{\eta}$. 

\begin{lmma}
\label{circlerotation}
Let $\alpha \in (0,1)$ be irrational and let $\eta \in (0,1)$ be given. 
\begin{enumerate}
\item[(1)] Suppose $x,y \in \bbt$ are distinct points. Then, there exists $m \in \bbz$ such that $h_{\eta}(R_{\alpha}^m(x)) \neq h_{\eta}(R_{\alpha}^m(y))$. 
\item[(2)] Suppose $\eta<\min\{\alpha,1-\alpha\}$ and let $x,y \in \bbt$ be distinct points. Then, there exists $m \in \bbz$ such that $\phi_{\eta,\alpha}(R_{\alpha}^{m}(x))\neq \phi_{\eta,\alpha}(R_{\alpha}^m(y))$. 
\end{enumerate}

\end{lmma}
\textit{Proof.} We denote the map $h_\eta:[0,1) \to \bbr$ simply by $h$.  Denote the periodic extension of $h$ to $\bbr$ by $\widetilde{h}$.  Suppose, for $x,y \in \bbt$,  $h(R_\alpha^{m}(x))=h(R_\alpha^{m}(y))$ for every $m \in \bbz$.  Then,  \[\widetilde{h}(x+m\alpha+n)=\widetilde{h}(y+m\alpha+n)\]
for all $m,n \in \bbz$.  This means that for all $m,n \in \bbz$, $(T_x\widetilde{h})(m\alpha+n)=(T_y\widetilde{h})(m\alpha+n)$, where $T_x\widetilde{h}, T_y\widetilde{h}$ denote the translations of $\widetilde{h}$ by $x,y$ respectively. 

Since the set $D=\{m\alpha+n: m,n \in \bbz\}$ is dense in $\bbr$, $T_x\widetilde{h}$, $ T_y\widetilde{h}$ are right continuous and they agree on $D$, we can conclude that they agree on $\bbr$.  Then, $T_z\widetilde{h}=\widetilde{h}$ on $\bbr$, where $z=x-y$. This clearly implies that   $z$ is an integer. This proves $(1)$. 

Let $\eta<\min\{\alpha,1-\alpha\}$ be given. Denote $\phi_{\eta,\alpha}$ by $\phi$. Note that for $x \in \bbt$, 

\begin{equation}
\label{functiontypeiii}
\phi(R_{\alpha}^{-1}x):=\begin{cases}
 1  & \mbox{ if
} x  \in [0,\eta),\cr
0 & \mbox{~if~} x \in [\eta,\alpha), \cr
    -1 &  \mbox{ if } x \in [\alpha,\alpha+\eta), \cr
    0 & \mbox{~if~} x \in [\alpha+\eta,1).
             \end{cases}
\end{equation}
Let $\chi:\{1,0,-1\} \to \{0,1\}$ be defined by $\chi(1)=1$ and $\chi(j)=0$ if $j \in \{0,-1\}$. Then, 
\[
\chi(\phi(R_{\alpha}^{-1}x))=h_{\eta}(x).\]
The proof of $(2)$ follows from the above equality and $(1)$. This completes the proof. \hfill $\Box$

Let $\alpha \in (0,1)$ be irrational and let $\eta \in (0,1)$ be such that $\eta < \min\{\alpha,1-\alpha\}$. Denote $h_{\eta}$ by $h$ and $\phi_{\eta,\alpha}$ by $\phi$. Let $Y_{\alpha,\eta}:=\bbt \times \bbz$ and let $X_{\alpha,\eta}:=\bbt \times \bbn$. Define a $\bbz^{2}$-action on $Y_{\alpha,\eta}$ by 
\[
(x,t)+e_1:=(x+\alpha, \phi(x)+t+1);~\textrm{and~~} (x,t)+e_2:=(x,t+1).\]

Since $\phi \geq -1$, it follows that $X_{\alpha,\eta}+\bbn^2 \subset X_{\alpha,\eta}$. Also, it is routine to verify that $(Y_{\alpha,\eta},X_{\alpha,\eta})$ is a pure $(\bbz^2,\bbn^2)$-space. 

Let $\beta$ be an arbitrary positive real number. Define a probability measure $\nu_{\alpha,\eta}$ on $\bbt$ by 
\[
\nu_{\alpha,\eta}(E):=\frac{1}{||e^{\beta h}||_1}\int_{E}e^{\beta h(x)}dx.\]
In the above formula, $dx$ is the Haar measure on the circle. Clearly, 
\[
\frac{d(\nu_{\alpha,\eta} \circ R_\alpha)}{d\nu_{\alpha,\eta}}=e^{\beta \phi}.\]

Define a measure $\mu_{\alpha,\eta}$ on $Y_{\alpha,\eta}$ by 
\[
\mu_{\alpha,\eta}(E \times \{n\})=(1-e^{-\beta})e^{-\beta n}\nu_{\alpha,\eta}(E).\]
It is easily verifiable that $\mu_{\alpha,\eta}$ is an $e^{-\beta c}$-conformal measure on the $(\bbz^2,\bbn^2)$-space $(Y_{\alpha,\eta},X_{\alpha,\eta})$. 

\begin{ppsn}
\label{type II rational}
	With the foregoing notation, we have the following. 
	\begin{enumerate}
		\item The $\bbz^2$-action on $Y_{\alpha,\eta}$ is free and ergodic.
		\item The measure $\mu_{\alpha,\eta}$ is of type II. 
		\item The map
		\[
		Y_{\alpha} \ni {(x,t)} \to Q_{(x,t)} \in Y_u\]
		is injective.
		\item Suppose $\alpha_1,\alpha_2 \in (0,\frac{1}{2})$ are distinct irrationals and let $\eta_1,\eta_2 \in (0,1)$ be such that for $i=1,2$, $\eta_i < \min\{\alpha_i,1-\alpha_i\}$.  Then, the $(\bbz^2,\bbn^2)$-spaces $(Y_{\alpha_1,\eta_1}, X_{\alpha_1,\eta_1},\mu_{\alpha_1,\eta_1})$ and 
		$(Y_{\alpha_2,\eta_2}, X_{\alpha_2,\eta_2},\mu_{\alpha_2,\eta_2})$ are not isomorphic. In particular, there are uncountably many type II extremal $\beta$-KMS states for $\theta=1$.
			\end{enumerate}  
\end{ppsn}
\textit{Proof.} Note that $\mu_{\alpha,\eta}$ is absolutely continuous w.r.t. to $dxdn$ which is a $\bbz^2$-invariant measure. The proof of $(1)$ and $(3)$ are very similar to the proof of Prop. \ref{adding machine}. The proof that the map 
\[
Y_{\alpha,\eta} \ni (x,t) \to Q_{(x,t)} \in Y_u\]
is injective is similar to Prop. \ref{adding machine}, where instead of Lemma \ref{separation}, we use Part $(2)$ of Lemma \ref{circlerotation}. 

We now prove $(4)$. Let $\alpha_1,\alpha_2 \in (0,\frac{1}{2})$ be irrationals and let $\eta_1,\eta_2 \in (0,1)$ be such that for $i=1,2$, $\eta_i < \min\{\alpha_i,1-\alpha_i\}$. Suppose $(Y_{\alpha_1,\eta_1},X_{\alpha_1,\eta_1},\mu_{\alpha_1,\eta_1})$
 and $(Y_{\alpha_2,\eta_2},X_{\alpha_2,\eta_2},\mu_{\alpha_2,\eta_2})$ are isomorphic. 
Then there exist $\bbz^2$-invariant null sets $N_{\alpha_i} \subset Y_{\alpha_i,\eta_i}$ for $i=1,2$ and an invertible measurable map $S:Y_{\alpha_1,\eta_1}\backslash N_{\alpha_1} \to Y_{\alpha_2,\eta_2} \backslash N_{\alpha_2}$ such that 
\begin{enumerate}
	\item[(i)] the map $S$ is $\bbz^2$-equivariant, $S(X_{\alpha_1,\eta_1}\backslash N_{\alpha_1})=X_{\alpha_2,\eta_2}\backslash N_{\alpha_2}$, and
	\item[(ii)] for every Borel subset $E \subset Y_{\alpha_2,\eta_2} \backslash N_{\alpha_2}$,  $\mu_{\alpha_2,\eta_2}(E)=\mu_{\alpha_1,\eta_1}(S^{-1}(E))$.
\end{enumerate}
Since $\mu_{\alpha_i,\eta_i}$ is absolutely continuous w.r.t. $dxdn$, it follows that $S$ preserves the measure $dxdn$. 
Let $S(x,t)=(S_1(x,t), S_2(x,t))$, where $S_1$ and $S_2$ are the coordinate functions. Note that  $e_2$ acts by  translation by 1 on the second coordinate. Hence, for $i=1,2$, $N_{\alpha_i}$ must be of the form $N_{\alpha_i}=U_{\alpha_i} \times \bbz$ for some $R_{\alpha_i}$-invariant subset $U_{\alpha_i}\subseteq \bbt$.
Since the map $S$ is $\bbz^2$-equivariant, from the action of $e_2$, we notice that, 
\begin{eqnarray}\label{v2action}
	S_1(x,t)=S_1(x,t+1)\quad \text{and}\quad S_2(x,t+1)=S_2(x,t)+1.
\end{eqnarray}
 In particular, the function $(x,t) \to S_1(x,t)$ does not depend on $t$.  As $S$ maps $X_{\alpha_1,\eta_1}\backslash N_{\alpha_1}$ onto $X_{\alpha_2,\eta_2}\backslash N_{\alpha_2}$ and $X_{\alpha_i,\eta_i}\backslash(X_{\alpha_i,\eta_i}+e_2)=\bbt \times \{0\}$,  it follows that $S_2(x,0)=0$ for $x \in \bbt \backslash U_{\alpha_1}$.   Eq. \ref{v2action} implies that $S_2(x,t)=t$ for $(x,t) \in \bbt\backslash U_{\alpha_1} \times \bbz$. 
 
Now, from the action of $e_1$, and by Equation \ref{v2action}, we can conclude that 
\[R_{\alpha_2} S_1(x,0)=S_1(R_{\alpha_1}x,\phi_{\eta_1,\alpha_1}(x)+1)=S_{1}(R_{\alpha_1}x,0).\]

Define a map $\widetilde{S}:\bbt \backslash U_{\alpha_1} \to \bbt \backslash U_{\alpha_2}$ by $\widetilde{S}(x)=S_1(x,0)$. Then, $\widetilde{S}$ is an isomorphism between $(\bbt,dx,R_{\alpha_1})$ and $(\bbt,dx,R_{\alpha_2})$. This implies that $\alpha_1=\alpha_2$. This completes the proof. \hfill $\Box$

\textbf{Type III examples:} The uncountably many type III examples that we construct are directly based on the works of Nakada (\cite{Nakada} and \cite{Nakada1}) on cylinder flows. We follow the exposition given in \cite{Aaronson_Nakada} and  recall a few facts from \cite{Aaronson_Nakada} concerning the ergodicity of cylinder flows. Let $\beta>0$ be fixed. For an irrational $\alpha \in (0,1)$, let $R_\alpha$ be the rotation on $\mathbb{T}$ by angle $\alpha$. We use additive notation. Thus, $\bbt=[0,1)$ and $R_\alpha(x)=x+\alpha \mod 1$. For $\gamma>0$, let $F_{\gamma}:[0,1) \to \bbc$ be defined by 
\begin{equation}
\label{functiontypeiii}
F_{\gamma}(x):=\begin{cases}
 1  & \mbox{ if
} x  \in [0,\frac{\gamma}{\gamma+1}),\cr
    -\gamma &  \mbox{ if } x \in [\frac{\gamma}{\gamma+1},1).
         \end{cases}
\end{equation}

\begin{rmrk}[Nakada]
\label{Nakada}
Thanks to Prop. 1.1, Thm. 1.4 and Thm. 1.6 of \cite{Aaronson_Nakada}, we have the following.
\begin{enumerate}
\item[(a)] There exists a unique non-atomic probability measure $m:=m_{\alpha,\beta,\gamma}$ on $\bbt$ such that $\frac{d(m\circ R_\alpha)}{dm}=e^{-\beta F_\gamma}$. Moreover, $(R_\alpha,m)$ is ergodic. 
\item[(b)] The measure $m_{\alpha,\beta,1}$ is of type III.
\item[(c)] Suppose $\alpha$ has bounded  partial quotients and $\frac{\gamma}{\gamma+1} \notin \mathbb{Q}+\mathbb{Q}\alpha$. Then, $m_{\alpha,\beta,\gamma}$ is of type III. 
\end{enumerate}
We also need the description of the spectrum of $R_\alpha$ given in \cite{Nakada1} (Page 476, Paragraph 1). Let $\beta, \gamma>0$ be given. Define 
\[
Spec(R_\alpha):=\{\lambda \in \bbt: \textrm{~There exists $\xi \in L^{\infty}(\bbt,m_{\alpha,\beta,\gamma})$ such that $\xi \circ R_\alpha=\lambda \xi$}\}.\]
If $\alpha$ has bounded partial quotients, then $Spec(R_\alpha)=\{e^{2\pi i n\alpha}:n \in \bbz\}$. In particular, if $\alpha_1,\alpha_2 \in (0,\frac{1}{2})$ are distinct irrationals having bounded partial quotients, then the dynamical systems 
$(\bbt,m_{\alpha_1,\beta,\gamma}, R_{\alpha_1})$ and $(\bbt,m_{\alpha_2,\beta,\gamma},R_{\alpha_2})$ are not metrically isomorphic. 
\end{rmrk}

\textbf{The case $\theta=1$:}
Let $\alpha \in (0,1)$ be an irrational and let $R_\alpha$ denote the rotation on $\bbt$ by angle $\alpha$. 
Take $\gamma=1$ in Eq. \ref{functiontypeiii} and consider $F_1$ given by,
\begin{equation*}
F_1(x)=\begin{cases}
1 & \mbox{ if~~} x \in [0,\frac{1}{2}),\\
-1 &  \mbox{ if } x \in [\frac{1}{2},1).
\end{cases}
\end{equation*} 
Then by $(a)$ and $(b)$ of Remark \ref{Nakada}, there exists a unique non-atomic type III probability measure $m:=m_{\alpha,\beta,1}$ such that $\frac{d(m\circ R_\alpha)}{dm}=e^{-\beta F_1}$. 

Recall that the homomorphism $c:G \to \bbr$ is given by $c(v_1)=1$ and $c(v_2)=2$. Define $Y_{\alpha}=\bbt \times \bbz$ and set $X_{\alpha}=\bbt \times \{0,1,2,\cdots\}$. Let $\mu_\alpha$ be the measure  on $Y_{\alpha}$ defined by,
$$\mu_\alpha(E\times \{n\})=(1-e^{-2\beta})e^{-2n\beta }m(E),$$
for a measurable subset $E\subset \bbt$.
Clearly, $\mu_{\alpha}(X_\alpha)=1$.
Define a $\bbz^2$-action on $Y_{\alpha}$ by 
\[
(x,t)+v_1:=(\tau_\alpha (x), 1_{[\frac{1}{2},1)}(x)+t);~\textrm{and~~} (x,t)+v_2:=(x,t+1).\]
It is easy to see that $(Y_{\alpha},X_{\alpha})$ is a pure $(\bbz^2,\bbn^2)$-space and $\mu_\alpha$ is an $e^{-\beta c}$-conformal measure on $Y_{\alpha}$.

\hfill $\Box$
\begin{ppsn}
\label{type III rational}
	With the foregoing notation, we have the following. 
	\begin{enumerate}
		\item The $\bbz^2$-action on $Y_\alpha$ is free and ergodic.
		\item The measure $\mu_\alpha$ is of type III. 
		\item The map
		\[
		Y_{\alpha} \ni {(x,t)} \to Q_{(x,t)} \in Y_u\]
		is injective.
		\item If $\alpha_1\neq \alpha_2$ are distinct irrationals in $(0,\frac{1}{2})$ having bounded partial quotients, then the $(\bbz^2,\bbn^2)$-spaces $(Y_{\alpha_1}, X_{\alpha_1},\mu_{\alpha_1})$ and 
		$(Y_{\alpha_2}, X_{\alpha_2},\mu_{\alpha_2})$ are not isomorphic. In particular, there are uncountably many type III extremal $\beta$-KMS states for $\theta=1$.
			\end{enumerate}  
\end{ppsn}
\textit{Proof.} With Remark \ref{Nakada} in hand,  the proof of $(1)$, $(2)$ and $(3)$ are exactly similar to the proof of Prop. \ref{adding machine} and hence we omit the proof. For example, the proof for the injectivity of the map  
\[
Y_{\alpha} \ni (x,t) \to Q_{(x,t)} \in Y_u\]
 is similar to that of Prop. \ref{adding machine}, where we work in the coordinate system determined by $\{v_1,v_2\}$,  the function $\phi$ is replaced with  $1_{[\frac{1}{2},1)}$ and in place of Lemma \ref{separation}, we appeal to  Lemma \ref{circlerotation} with $\eta=\frac{1}{2}$. With the description of the spectrum alluded to in Remark \ref{Nakada}, the proof of $(4)$ is very similar to the proof of Part $(4)$ of Lemma \ref{type II rational}. We leave the details to the reader. \hfill $\Box$

\textbf{The case of an irrational $\theta$:}
The construction of type III examples for an irrational $\theta$ is similar to the case $\theta=1$. Recall that the homomorphism $c:G \to \bbr$ is given by $c(e_1)=1$ and $c(e_2)=\theta$. Take $\gamma=\theta$ in \ref{functiontypeiii}, and let $F_{\theta}$ be the function defined by 
\begin{equation*}
F_{\theta}(x)=\begin{cases}
1 & \mbox{ if~~} x \in [0,\frac{\theta}{\theta+1}),\\
-\theta &  \mbox{ if } x \in [\frac{\theta}{\theta+1},1).
\end{cases}
\end{equation*} 
 Choose $\alpha \in [0,1)\backslash\mathbb{Q}$ such that $\alpha$ has  bounded partial quotients and $\frac{\theta}{\theta+1} \notin \mathbb{Q}+\mathbb{Q}\alpha$. Then, by $(a)$ and $(c)$ of Remark \ref{Nakada}, there exists a unique non-atomic type III probability measure $m:=m_{\alpha,\beta,\theta}$ on $\bbt$ such that $\frac{d(m\circ R_\alpha)}{dm}=e^{-\beta F_\theta}$.  
 
 Define $Y_{\alpha}=\bbt \times \bbz$ and set $X_{\alpha}=\bbt \times \{0,1,2,\cdots\}$. Let $\mu_\alpha$ be the measure on $Y_{\alpha}$ given by,
$$\mu_\alpha(E\times \{n\})=(1-e^{-\beta(1+\theta)})e^{-\beta n (\theta+1)}m(E)$$
for a measurable subset $E\subset \bbt$. 
Clearly, $\mu_{\alpha}(X_{\alpha})=1$. 
Define a $\bbz^2$-action on $Y_{\alpha}$ by 
\[
(x,t)+v_1:=(R_\alpha (x), 1_{[\frac{\theta}{\theta+1},1)}(x)+t);~\textrm{and~~} (x,t)+v_2:=(x,t+1).\]
It is easy to see that $(Y_{\alpha},X_{\alpha})$ is a pure $(\bbz^2,\bbn^2)$-space and $\mu_\alpha$ is an $e^{-\beta c}$-conformal measure on $Y_\alpha$.

If we make use of Remark \ref{Nakada}, the proof of the following proposition is similar to that of Prop. \ref{type III rational}, and hence the proof is omitted. 
\begin{ppsn}
With the foregoing notation, we have the following. 
\begin{enumerate}
		\item The $\bbz^2$-action on $Y_\alpha$ is free and ergodic.
		\item The measure $\mu_\alpha$ is of type III. 
		\item The map
		\[
		Y_{\alpha} \ni {(x,t)} \to Q_{(x,t)} \in Y_u\]
		is injective.
		\item Suppose $\alpha_1\neq \alpha_2$ are distinct irrationals in $(0,\frac{1}{2})$ having bounded partial quotients such that $\frac{\theta}{1+\theta} \notin \mathbb{Q}+\mathbb{Q}\alpha_i$ for $i=1,2$.  Then, the $(\bbz^2,\bbn^2)$-spaces $(Y_{\alpha_1}, X_{\alpha_1},\mu_{\alpha_1})$ and $(Y_{\alpha_2}, X_{\alpha_2},\mu_{\alpha_2})$ are not isomorphic. In particular, there are uncountably many type III extremal $\beta$-KMS states when $\theta$ is irrrational. 
			\end{enumerate}  
\end{ppsn}

The type II and type III examples for other rational values of $\theta$ can be constructed from the case $\theta=1$ as follows. Let $\theta:=\frac{p}{q}$ be a positive rational. Assume that $\gcd(p,q)=1$. Choose a matrix $\begin{bmatrix}
x & y \\
z & w\\
\end{bmatrix} \in SL_2(\bbz)$ such that $x,y,z,w \geq 0$ and 
\[
x+z=q;~y+w=p.\]
To see that this is possible, observe that if $q=1$, then the matrix $\begin{bmatrix}
1 & p-1 \\
0 & 1
\end{bmatrix}$ will do. Suppose $q \geq 2$. Choose $x \in \{1,2,\cdots,q-1\}$ such that $xp \equiv 1 \mod q$. Write $xp=yq+1$ with $y \in \bbz$. Set $z:=q-x$ and $w:=p-y$. Then, $x,y,z,w \geq 0$, $\begin{bmatrix}
x & y \\
z & w 
\end{bmatrix} \in SL_2(\bbz)$ and 
\[
x+z=q;~y+w=p.\]

Let us fix notation. Let $\phi:\bbz^2 \to \bbz^2$ be the isomorphism that corresponds to the matrix $\begin{bmatrix}
x & y \\
z & w
\end{bmatrix}$. Note that $\phi(\bbn^2) \subset \bbn^2$. We denote the homomorphism $\bbz^2 \to \bbr$ that sends $e_1 \to 1$ and $e_2 \to \theta$ by $c_{\theta}$. If $\theta=1$, we denote $c_\theta$ simply by $c$. It is clear that $\frac{1}{q}( c \circ \phi)=c_\theta$.

Fix $\beta>0$ and let $t \in \{II,III\}$. Let $(Y,X)$ be a pure $(\bbz^2,\bbn^2)$-space together with an $e^{-\frac{\beta c}{q}}$-conformal measure $m$ on $Y$ satisfying $(A1)$ and $(A2)$ of Remark \ref{defining example}. Define 
\[
Y^{\phi}:=Y;~~X^{\phi}:=X;~~m^{\phi}:=m.\]
Define a new $\bbz^2$-action on $Y^{\phi}$ as follows:
\[
y \oplus (m,n):=y+\phi(m,n).\]
With this new action, $(Y^{\phi},X^{\phi},m^{\phi})$ is a pure $(\bbz^2,\bbn^2)$-space and $m^{\phi}$ is $e^{-\beta c_\theta}$-conformal. Since $\phi$ is an automorphism of $\bbz^2$, it follows that $(Y^{\phi},m^{\phi},\bbz^2)$ inherits the freeness, ergodicity and the type from that of $(Y,m,\bbz^2)$. Note that for $y \in Y^{\phi}$, 
\[
Q^{\phi}_y:=\{(m,n) \in \bbz^2: y \ominus (m,n) \in X^{\phi}\}=\phi^{-1}(Q_y).\]
Since the map $Y \ni y \to Q_y \in Y_u$ is injective, it follows that $Y^{\phi} \ni y \to Q_{y}^{\phi} \in Y_u$ is injective. Thus, the $(\bbz^2,\bbn^2)$-space $(Y^{\phi},X^{\phi})$ together with the $e^{-\beta c_\theta}$-conformal measure satisfies $(A1)$ and $(A2)$ of Remark \ref{defining example}. 

Moreover, since $\phi$ is an automorphism, it is clear that the above construction, i.e. 
\[
(Y,X,m) \to (Y^{\phi},X^{\phi},m^{\phi})\]
maps metrically non-isomorphic spaces to metrically non-isomorphic spaces. 

As we have already proved the existence of a continuum of metrically non-isomorphic $(\bbz^2,\bbn^2)$-spaces satisfying $(A1)$ and $(A2)$ for $\theta=1$ and for every $\beta>0$, we can conclude that that there is a continuum of metrically non-isomorphic $(\bbz^2,\bbn^2)$-spaces that satisfy $(A1)$ and $(A2)$ for every rational $\theta>0$ and for every $\beta>0$.  We have now proved the following theorem. 

\begin{thm}
Suppose $\theta>0$
Let $\sigma:=\{\sigma_t\}_{t \in \bbr}$ be the $1$-parameter group of automorphisms on $C_{c}^{*}(\bbn^2)$ given by 
\[
\sigma_t(v_{(m,n)})=e^{i(m+n\theta)t}v_{(m,n)}.\]
Then, for every $\beta>0$, there is a continuum of extremal $\beta$-KMS states on $C_{c}^{*}(\bbn^2)$ for $\sigma$ of both type II and type III. 
\end{thm}

We end our paper by posing a problem. Constructing type III product measures for the Bernoulli shift on two symbols have always been a subject of interest. The first such example was due to Hamachi (\cite{Hamachi}). We refer the reader to \cite{Vaes_Kosloff} and the references therein for more recent developments. In the recent years, the focus is on constructing product measures of various Krieger types. Given this interest, the following question, we believe, is worth investigating. 

\textbf{Question:} Let $\theta>0$ be irrational and suppose $\beta>0$. Let $\chi:\{0,1\}^{\bbz} \to \bbr$ be the potential defined by \begin{equation*}
\chi(x)=\begin{cases}
1 & \mbox{ if~~} x_{0}=0,\\
-\theta &  \mbox{ if } x_{0}=1.
\end{cases}
\end{equation*} 
\begin{enumerate}
\item[(1)] Does there exist an ergodic probability measure  for the Bernoulli shift of type $II_\infty$ that is $e^{-\beta \chi}$-conformal ?
\item[(2)] Let $\lambda \in [0,1)$ be given. Does there  exist an ergodic probability measure of type $III_\lambda$ that is $e^{-\beta \chi}$-conformal ?
\end{enumerate} 
Nakada's work shows that type $III_1$ is possible. Note that the type II examples that we have constructed in this paper are of type $II_1$.

\bibliography{references}
 \bibliographystyle{amsplain}
 
 \vspace{2.5 mm}
 \noindent
  {\sc Anbu Arjunan}
 (\texttt{aanbu100@gmail.com})\\
 {\footnotesize Institute of Mathematical Sciences,  \\A CI of Homi Bhabha National Institute, \\
     CIT Campus, Taramani, Chennai, 600113,\\ Tamilnadu, India.}\\
 
 \noindent 
  {\sc Sruthymurali}
  (\texttt{sruthy92smk@gmail.com})\\
  {\footnotesize Chennai Mathematical Institute, \\
  Siruseri IT Park, 
    Siruseri, Chennai, 603103, \\Tamilnadu, India.}\\
 
 \noindent
{\sc S. Sundar}
(\texttt{sundarsobers@gmail.com})\\
         {\footnotesize  Institute of Mathematical Sciences, \\
         A CI of Homi Bhabha National Institute, \\
CIT Campus, Taramani, Chennai, 600113, \\Tamilnadu, INDIA.}\\

\end{document}